%% file: paper_kdro.tex
\documentclass[twoside]{article}

\usepackage[accepted]{aistats2021}
\usepackage[utf8]{inputenc} 
\usepackage[T1]{fontenc}    
\usepackage{hyperref}       
\usepackage{url}            
\usepackage{booktabs}       
\usepackage{amsfonts}       
\usepackage{nicefrac}       
\usepackage{microtype}      
\usepackage{subfig}

\newcommand{\ourtitle}{Kernel Distributionally Robust Optimization
}

\usepackage[pdftex]{graphicx}

\usepackage[authoryear]{natbib}
\setcitestyle{authoryear,open={(},close={)}}
\renewcommand{\cite}{\citep}

\input{draft/preamble.tex}

\begin{document}
\runningtitle{Kernel Distributionally Robust Optimization}
\runningauthor{Jia-Jie Zhu, Wittawat Jitkrittum, Moritz Diehl, Bernhard Sch\"olkopf}
\twocolumn[

\aistatstitle{\ourtitle}

\aistatsauthor{ Jia-Jie Zhu \And Wittawat Jitkrittum }

\aistatsaddress{ 
  Empirical Inference Department\\
  Max Planck Institute for Intelligent Systems\\
   T\"ubingen, Germany\\
   \texttt{jia-jie.zhu@tuebingen.mpg.de}\\
   \And  
   Empirical Inference Department\\
   Max Planck Institute for Intelligent Systems\\
    T\"ubingen, Germany\\
    Currently at Google Research, NYC, USA\\
    \texttt{wittawatj@gmail.com}
    }    
  
\aistatsauthor{ Moritz Diehl \And
    Bernhard Sch\"olkopf}
    \aistatsaddress{
   Department of Microsystems Engineering\\
   \& Department of Mathematics\\
   University of Freiburg\\
   Freiburg, Germany\\
   \texttt{moritz.diehl@imtek.uni-freiburg.de}
   \And 
Empirical Inference Department\\
Max Planck Institute for Intelligent Systems\\
 T\"ubingen, Germany\\
 \texttt{bernhard.schoelkopf@tuebingen.mpg.de}
} ]

\begin{abstract}
We propose \emph{kernel distributionally robust optimization} (\kdro) using insights from the robust optimization theory and functional analysis. Our
method uses reproducing kernel Hilbert spaces (RKHS) to construct a wide range of convex ambiguity sets, which can be generalized to sets based on integral probability metrics and finite-order moment bounds. This perspective unifies multiple existing robust and stochastic optimization methods. We prove a theorem that generalizes the classical duality in the mathematical problem of moments. Enabled by this theorem, we reformulate the maximization with respect to measures in DRO into the dual program that searches for RKHS functions. Using universal RKHSs, the theorem applies to a broad class of loss functions, lifting common limitations such as polynomial losses and knowledge of the Lipschitz constant. We then establish a connection between DRO and stochastic optimization with expectation constraints. Finally, we propose practical algorithms based on both batch convex solvers and stochastic functional gradient, which apply to general optimization and machine learning tasks.
\end{abstract}

\section{INTRODUCTION}
\input{draft/intro.tex}
\section{BACKGROUND}

\paragraph{Notation.}
$\mathcal X \subset \mathbb{R}^d$ denotes the input domain, which is 
assumed to be compact unless otherwise specified.
$\mathcal P:=\mathcal P(\mathcal X) $ denotes the set of all Borel probability measures on $\mathcal X$. 
We use ${\hat P}$ to denote the empirical distribution ${\hat P}=\sum_{i=1}^N{\frac1N}\delta_{\xi_i}$, where $\delta$ is a Dirac measure and $\{\xi_i\}_{i=1}^N$ are data samples.
We refer to the function $\indc (x): = 0$ if $x\in \mathcal C$, $\infty$ if $x\notin \mathcal C$, as the indicator function. 
$\suppc (f):=\sup_{\mu\in\mathcal C} \hip{f}{\mu} $ is the support function of $\mathcal C$.
$S_{N}$ denotes the $N$-dimensional simplex.
$\mathrm{ri}(\cdot)$ denotes the relative interior of a set.
A function $f$ is upper semicontinuous on \domain if $\limsup_{x\to x_0} f(x) \leq f(x_0), \forall x_0\in\domain$; it is proper if it is not identically $-\infty$.
When there is no ambiguity, we simplify the loss function notation $l(\var,\cdot)$ by using $l$ to indicate that results hold for \var point-wise.
\input{draft/background.tex}

\section{THEORY}
\label{sec:theory}
\input{draft/conic.tex}

\section{COMPUTATION}
\label{sec:practical}
\input{draft/practical.Rmd}
\input{draft/aistats/sgd_kdro.tex}
\input{draft/aistats/exp_fig_all.tex}
\section{NUMERICAL STUDIES}
\label{sec:numerical}
This section showcases the applicability of \kdro (and hence \sfg) and discusses the robustness-optimality trade-off.
Our purpose is not to benchmark state-of-art performances or to demonstrate the superiority of a specific algorithm. Indeed, we believe both RO and DRO are elegant theoretical frameworks that have their specific use cases.
We note that our theory can be applied to a broader scope of applications than the examples here, such as stochastic optimal control.
See the appendix for more experimental results.
The code is available at \url{https://github.com/jj-zhu/kdro}.
\paragraph{Distributionally robust solution to uncertain least squares.}
\input{draft/exp_1.tex}
\input{draft/aistats/exp_adv.tex}

\section{OTHER RELATED WORK AND DISCUSSION}
\label{sec:lit}
\input{draft/literature.tex}
\input{draft/discuss.Rmd}
\subsubsection*{Acknowledgements}
  We thank Daniel Kuhn for the helpful discussion during a workshop at IPAM, UCLA.
  We also thank Yassine Nemmour and Simon Buchholz for sending us their feedback on the paper draft.
  During part of this project, Jia-Jie Zhu was supported by the European Union’s Horizon 2020 research and innovation programme under the Marie Skłodowska-Curie grant agreement No 798321.
  Moritz Diehl would like to acknowledge the funding support from the
  DFG via Project DI 905/3-1 (on robust model predictive control)
\bibliography{paper.bib}

\clearpage
\newpage
\appendix
\onecolumn

\begin{center}
{\LARGE{}Appendix: \ourtitle{}}{\LARGE\par}
\par\end{center}


\section{PROOFS OF THEORETICAL RESULTS}
\label{sec:proof}
\input{draft/proof.Rmd}

\section{COMPUTATIONAL FORMULATIONS}
\input{draft/aistats/pluggin.md}
\subsection{Random features}
\input{draft/randfeat.md}

\section{FURTHER NUMERICAL EXPERIMENT RESULTS}
\label{sec:exp_more}
\input{draft/certify.Rmd}
\input{draft/pf/gap.Rmd}

\section{SUPPORTING LEMMAS}
\label{sec:lemmas}
\input{draft/compact.Rmd}

\end{document}

%% file: draft/preamble.tex
\usepackage{amsthm}
\usepackage{amsmath}
\usepackage{amsfonts}
\usepackage{bbm}
\usepackage{xspace}
\usepackage{tikz}
\usepackage{enumitem} 
\usepackage{graphbox} 
\usepackage{soul}
\usepackage{wrapfig}
\usepackage[export]{adjustbox} 
\usepackage{algorithm}
\usepackage[noend]{algorithmic}
\usepackage{float}
\usetikzlibrary{decorations.pathreplacing}

\newif\ifcomment
\commenttrue

\newcommand{\kdro}{Kernel DRO\xspace} 
\newcommand{\sfg}{SFG-DRO\xspace}
\newcommand{\var}{\ensuremath{\theta}\xspace} 
\newcommand{\pfeas}{\ensuremath{\mathcal{K}_{\mathcal C}}\xspace} 

\newcommand{\hnorm}[1]{\|{#1}\|_{\mathcal H}}
\newcommand{\hip}[2]{\langle {#1, #2} \rangle _{\mathcal H}}
\newcommand{\pip}[2]{ {\int #1\ d #2}} 
\newcommand{\expected}[2]{{\mathbb E}_{#1}\ #2}
\newcommand{\rkhs}{\ensuremath{\mathcal H}\xspace}
\newcommand{\kme}[2]{\sum_{i=1}^N{#1} \phi({#2})}
\newcommand{\mnz}{\min} 
\newcommand{\mxz}{\max}
\newcommand{\indc}{\ensuremath{\delta_{\mathcal C}}} 
\newcommand{\suppc}{\ensuremath{\delta^*_{\mathcal C}}} 
\newcommand{\semiinfcons}{\ensuremath{l(\var, \xi)\leq  f_0 + f(\xi)   ,\ \forall \xi \in\mathcal X}} 

\newcommand{\rset}{\mathbb R}  
\newcommand{\co}[1]{\mathrm{co}(#1)} 
\newcommand{\domain}{\ensuremath{\mathcal X}\xspace}
\newcommand{\ipm}{\ensuremath{d_{\mathcal F}}\xspace} 
\newcommand{\mmap}{\ensuremath{\mathcal T}\xspace}

\DeclareMathOperator*{\sjt}{\mathrm{subject}\ \mathrm{to}}

\newtheorem{prop}{Proposition}[section]
\newtheorem{theorem}{Theorem}[section]
\newtheorem{lemma}[prop]{Lemma}
\newtheorem{assume}[prop]{Assumption}

\theoremstyle{definition}
\newtheorem{ex}{Example}[section]

\newtheorem{corol}{Corollary}[prop]

\newtheorem*{remark}{Remark}

\newenvironment{opt} 
    {
        \newcommand{\MIN}[2]{\min_{##1}& & & ##2\\}
        
        \newcommand{\SJT}[1]{\sjt & & &##1\\} 
        \begin{aligned}
    }
    {
        \end{aligned}
    }
\newenvironment{opt1} 
{
    \newcommand{\MINST}[3]{
        \begin{aligned}
          &\min_{##1}&&  ##2\\
          &\sjt&&\ ##3
          \end{aligned}
      }
    \newcommand{\MAXST}[3]{\max_{##1}  ##2\quad\sjt\ ##3}
    \newcommand{\SUPST}[3]{\sup_{##1}  ##2\quad\sjt\ ##3}
    
}
{
}

\makeatletter
\newcommand\footnoteref[1]{\protected@xdef\@thefnmark{\ref{#1}}\@footnotemark}
\makeatother

\usepackage{xcolor}

\ifcomment

\newcommand{\wjinline}[1]{
    [\textcolor{green!30!black}{\textbf{WJ:}} \textcolor{green!50!black}{\footnotesize\textbf{#1}}]
}
\newcommand{\jzinline}[1]{
    [\textcolor{blue!60!black}{\textbf{JZ:}} \textcolor{blue!60!black}{\footnotesize\textbf{#1}}]
}
\newcommand{\wjsay}[1]{
    \marginpar{
    [\textcolor{orange!60!black}{\textbf{WJ:}} \textcolor{orange!60!black}{\footnotesize\textbf{#1}}]
    }
    }
\newcommand{\jz}[1]{\jzinline{#1}}
\else
\newcommand{\letsremove}[1]{}
\newcommand{\wjsay}[1]{}
\newcommand{\wjinline}[1]{}
\newcommand{\jzinline}[1]{}
\newcommand{\jz}[1]{}
\fi

%% file: draft/intro.tex
Imagine a hypothetical scenario in the illustrative figure
where we want to arrive at a destination while
avoiding unknown obstacles. 
A \emph{worst-case robust optimization} (RO)~\cite{ben-talRobustOptimization2009a} approach is then to avoid the entire unsafe area (left, blue).
Suppose we have historical locations of the obstacles (right, dots).
We may choose to avoid only the convex polytope that contains all the samples (pink).
This \emph{data-driven robust decision-making} idea improves efficiency while retaining robustness.
\begin{figure}[h!]
    \vspace{-0.2cm}
    \centering
    \includegraphics[width=0.8\columnwidth]{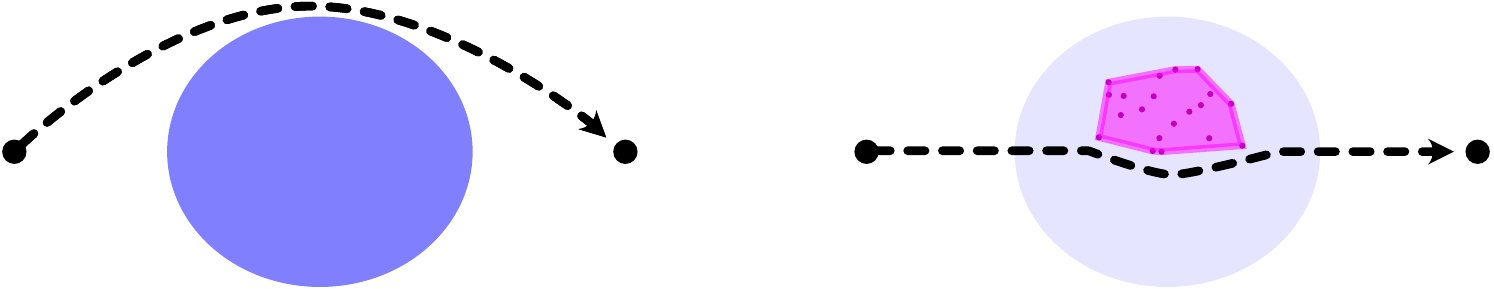}
    \label{fig:obstacles}
    \vspace{-0.45cm}
  \end{figure}

The concept of distributional ambiguity concerns the uncertainty of uncertainty --- the underlying probability measure is only partially known or subject to change. 
This idea is by no means a new one.
The classical moment problem concerns itself with estimating the worst-case risk expressed by $\max_{P\in\mathcal K}{\int l \, dP}$ where $l$ is some loss function.
The constraint $P\in\mathcal K$ describes the \emph{distribution ambiguity}, i.e., $P$ is only known to live within a subset $\mathcal K$ of probability measures.
The solution to the moment problem gives the risk under some worst-case distribution within $\mathcal K$.
To make decisions that will minimize this worst-case risk is the idea of \emph{distributionally robust optimization} (DRO)~\cite{delageDistributionallyRobustOptimization2010,scarfMinmaxSolutionInventory1958}.

Many of today's learning tasks suffer from various manifestations of distributional ambiguity --- e.g., covariate shift, adversarial attacks, simulation to reality transfer --- 
phenomena that are caused by the discrepancy between training and test distributions.
Kernel methods are known to possess robustness properties, e.g.,  
\cite{christmannConsistencyRobustnessKernelbased2007,xuRobustnessRegularizationSupport}.
However, this robustness only applies to kernelized models.
This paper extends the robustness of kernel methods
using the robust counterpart formulation techniques \cite{ben-talRobustOptimization2009a} as well as the principled conic duality theory \cite{shapiroDualityTheoryConic2001}. 
We term our approach \emph{kernel distributionally robust optimization} (\kdro), which can robustify general optimization solutions not limited to kernelized models.

The \emph{main contributions} of this paper are:
\begin{enumerate}[noitemsep,topsep=0pt]
  \item We rigorously prove the generalized duality theorem (Theorem~\ref{thm:dro_main_thm}) that reformulates general DRO into a convex dual problem searching for RKHS functions, lifting common limitations of DRO on the loss functions, such as the knowledge of Lipschitz constant.
  The theorem also constitutes a generalization of the duality results from the literature of mathematical problem of moments.
  \item We use RKHSs to construct a wide range of convex ambiguity sets (in Table~\ref{tbl:part},~\ref{tbl:sets}), including sets based on integral probability metrics (IPM) and finite-order moment bounds. This perspective unifies existing RO and DRO methods. 
  \item We propose computational algorithms based on both convex solvers and stochastic approximation, which can be applied to robustify general optimization and machine learning models not limited to kernelized or known-Lipschitz-constant ones.
  \item Finally, we establish an explicit connection between DRO and stochastic optimization with expectation constraints. This leads to a novel stochastic functional gradient DRO (\sfg) algorithm which can scale up to modern machine learning tasks.
\end{enumerate}


In addition, we give complete self-contained proofs in the appendix that shed light on the connection between RKHSs, conic duality, and DRO.
We also show that universal RKHSs are large enough for DRO from the perspective of functional analysis through concrete examples. 

%% file: draft/background.tex
\subsection{Robust and distributionally robust optimization}
\input{draft/ro_basics.Rmd}

\subsection{Reproducing kernel Hilbert spaces}
A symmetric function $k\colon \mathcal{X} \times \mathcal{X}\to \rset$ is called a positive definite kernel if $\sum_{i=1}^n \sum_{i=1}^n a_i a_j k(x_i, x_j)\ge 0$ for any $n \in \mathbb{N}$, $\{ x_i \}_{i=1}^n \subset \mathcal{X}$, and $\{a_i\}_{i=1}^n \subset \mathbb{R}$.
Given a positive definite kernel $k$, there exists a Hilbert space \rkhs and a feature map $\phi\colon
\mathcal{X} \to \rkhs$, for which $k(x,y) = \langle \phi(x), \phi(y)
\rangle_\rkhs$ defines an inner product on $\rkhs$, where $\rkhs$ is a space of real-valued functions on $\mathcal{X}$. The space $\rkhs$ is called a reproducing kernel Hilbert space (RKHS). It is equipped with the \emph{reproducing property}: 
$f(x) = \langle f, \phi(x) \rangle_\rkhs$ for any $f\in \rkhs, x \in
\mathcal{X}$.
By convention, we will denote the canonical feature map as $\phi(x):=k(x, \cdot)$.
Properties of the functions in $\rkhs$ are inherited from the properties of
$k$. For instance, if $k$ is continuous, then any $f \in \rkhs$ is continuous.
A continuous kernel $k$ on a compact metric space $\mathcal{X}$ is said to be
\emph{universal} if $\mathcal{H}$ is dense in $C(\mathcal{X})$ \cite[Section
4.5]{steinwartSupportVectorMachines2008}. A universal $\rkhs$ can thus be
considered a large RKHS since any continuous function can be approximated
arbitrarily well by a function in $\rkhs$. An example of a universal kernel is
the Gaussian kernel $k(x,y) = \exp\left( - \frac{\|x-y\|^2_2}{2\sigma^2}
\right)$ defined on $\mathcal{X}$ where
$\sigma>0$ is the bandwidth parameter.

RKHSs first gained widespread attention following the advent of the
kernelized support vector machine (SVM) for classification
problems \cite{CorVap1995,boserTrainingAlgorithmOptimal1992,scholkopfNewSupportVector2000}. More recently, the use of RKHSs has been extended to manipulating
and comparing probability distributions via kernel mean embedding
\cite{SmoGreSonSch2007}. Given a distribution $P$, and a (positive definite)
kernel $k$, the \emph{kernel mean embedding} of $P$ is defined as $\mu_P := \int
k(x, \cdot) \, dP$. If $\mathbb{E}_{x\sim P}[k(x,x)] < \infty$,
then $\mu_P \in \rkhs$ \cite[Section 1.2]{SmoGreSonSch2007}.
The reproducing property allows one to easily compute the expectation of any function $f
\in \rkhs$ since $\mathbb{E}_{x \sim P}[f(x)] = \langle f, \mu_P
\rangle_\rkhs$.
Embedding distributions into $\rkhs$ also allows one to
measure the distance between distributions in $\rkhs$.
If $k$ is universal, then the mean map $P \mapsto \mu_P$ is injective on $\mathcal{P}$
\cite{grettonKernelTwosampleTest2012}.
With a universal $\rkhs$, given two distributions $P,Q$, $\| \mu_P - \mu_Q \|_\rkhs$ defines a metric. This quantity is known as the maximum mean discrepancy (MMD) \cite{grettonKernelTwosampleTest2012}. 
With $\|f\|_\rkhs := \sqrt{\langle f, f \rangle_\rkhs }$ and the reproducing
property, it can be shown that $\| \mu_P - \mu_Q \|_\rkhs^2 = \mathbb{E}_{x,x'\sim
P}k(x,x') + \mathbb{E}_{y,y'
\sim Q} k(y, y') - 2\mathbb{E}_{x\sim P, y\sim Q} k(x,y)$, allowing the plug-in estimator to be used for estimating the MMD from empirical data.
The MMD is an instance of the class of integral probability metrics (IPMs),
and can equivalently be written as $\| \mu_P - \mu_Q \|_\rkhs = \sup_{
\|f\|_\rkhs \le 1} \int f \, d(P-Q)$, where the optimum $f^*$ is a witness function
\cite{grettonKernelTwosampleTest2012,sriperumbudurEmpiricalEstimationIntegral2012}.




%% file: draft/ro_basics.Rmd
Robust optimization (RO)~\cite{ben-talRobustOptimization2009a}
studies mathematical decision-making under uncertainty. It solves the min-max problem (omitting constraints)
$
\ \min_\var \sup_{\xi\in \domain}  l(\var,\xi),
$
where $l(\var, \xi)$ denotes a general loss function, $\var$ is the decision
variable, and $\xi$ is a variable representing the uncertainty.
Intuitively, RO makes the decision assuming an adversarial scenario, as reflected in taking supremum w.r.t. $\xi$. For this reason, it is often referred to as the \emph{worst-case RO}. 
Recently, RO has been applied to the setting of adversarially robust learning, e.g., in~\cite{madryDeepLearningModels2019,wongProvableDefensesAdversarial}, which we will visit in this paper.
In the optimization literature, a typical approach to solving RO is via
reformulating the min-max program using duality to obtain a single minimization problem. 
In contrast, distributionally robust optimization (DRO)
minimizes the expected loss assuming the worst-case distribution:
\begin{equation}
    \mnz_\var \underset{P\in\mathcal K}{\sup} \bigg\{ \int l(\var,\xi) \ dP(\xi) \bigg\},
    \label{eq:dro_intro}
\end{equation}
where $\mathcal K\subseteq \mathcal P$, called the \emph{ambiguity set}, is a subset of distributions, e.g., all distributions with the given mean and variance.
Compared with RO, 
DRO only robustifies the solution against a subset $\mathcal K$ of distributions on \domain and is, therefore, less conservative (since
$ \sup_{P\in\mathcal K} \{ \int l(\var,\xi) \ dP(\xi)  \}
\leq  \sup_{\xi\in \domain}  l(\var,\xi)
$).

The inner problem of \eqref{eq:dro_intro},
historically known as the \emph{problem of moments} traced back at least to Thomas Joannes Stieltjes, estimates the worst-case risk under uncertainty in distributions. 
The modern approaches, pioneered by the work of \cite{isiiSharpnessTchebychefftypeInequalities1962} (see also \cite{lasserreBoundsMeasuresSatisfying2002,
shapiroDualityTheoryConic2001,
bertsimasOptimalInequalitiesProbability2005,popescuSemidefiniteProgrammingApproach2005,vandenbergheGeneralizedChebyshevBounds2007,
vanparysGeneralizedGaussInequalities2016, zhuWorstCaseRiskQuantification2020}), typically seek a sharp upper bound via duality.
This duality, rigorously justified in \cite{shapiroDualityTheoryConic2001}, is different from that in the Euclidean space
because infinite-dimensional convex sets can become pathological. 
Using that methodology, we can reformulate DRO~\eqref{eq:dro_intro} into a single solvable minimization problem.

Existing DRO approaches can be grouped into three main categories by the type of ambiguity sets used.
DRO with (finite-order) moment constraints has been studied in \cite{delageDistributionallyRobustOptimization2010,scarfMinmaxSolutionInventory1958,zymlerDistributionallyRobustJoint2013}.
The authors of \cite{ben-talRobustSolutionsOptimization2013,iyengarRobustDynamicProgramming2005,nilimRobustControlMarkov2005,wangLikelihoodRobustOptimization2016a,duchiStatisticsRobustOptimization2016}
studied DRO using likelihood bounds as well as $\phi-$divergence.
Wasserstein-distance-based DRO has been studied by the authors of \cite{mohajerinesfahaniDatadrivenDistributionallyRobust2018,zhaoDatadrivenRiskaverseStochastic2018,gaoDistributionallyRobustStochastic2016,blanchetRobustWassersteinProfile2019}, and applied in a large body of literature.
Many existing approaches require either the assumptions such as quadratic loss functions or the knowledge of Lipschitz constant or RKHS norm of the loss $l$, which are often hard to obtain in practice; see \cite{virmauxLipschitzRegularityDeep2018,biettiKernelPerspectiveRegularizing2019}.

%% file: draft/conic.tex
We make the following assumption for the proof.
\begin{assume}
        \vspace{-0.04cm}
        \label{asm:compact}
        $l(\var, \cdot)$ is proper, upper semicontinuous.
        $\mathcal C$ is closed convex. $\mathrm{ri}(\pfeas)\neq \emptyset $. \
        \vspace{-0.2cm}
\end{assume}
This assumption is general in that it does not require the knowledge of the Lipschitz constant or the RKHS $l(\var, \cdot)$ lives in.
Generally speaking, the DRO problem~\eqref{eq:dro_intro} requires two essential elements: an appropriate ambiguity set that contains meaningful distributions and a sharp reformulation of the min-max problem. We first present the former in Section~\ref{sec:primal}, and then the latter in Section~\ref{sec:dual}.
Complete proofs of our theory are deferred to the appendix.
\subsection{Generalized primal formulation}
\label{sec:primal}
We now present
the primal formulation of kernel distributionally robust optimization (\kdro) as a generalization of existing DRO frameworks.
\begin{multline}
        (P):=\min_{\var} \underset{P,\mu}{\sup}\bigg\{{  \int l(\theta,\xi) \ dP(\xi) }\colon\\
        {\int{\phi }{\ dP} = \mu, P\in \mathcal P,\mu\in \mathcal C} \bigg\},
\label{eq:kdro}
\end{multline}
where $ \rkhs$ is an RKHS  whose feature map is $\phi$.
Both sides of the constraint $\int{\phi }{\ dP} = \mu$ are functions in $\rkhs$. Note $\mu$ can be viewed as a generalized moment vector, which is constrained to lie within the set $ \mathcal C\subseteq\rkhs$,
referred to 
as an (RKHS) ambiguity set.
Let us denote the set of all feasible distributions in \eqref{eq:kdro} as $ \pfeas = {\{P \colon \int \phi \ dP  = \mu, \mu\in\mathcal C, P\in \mathcal P\}}$, i.e., \pfeas is the usual ambiguity set.
Intuitively, the set $\mathcal C$ restricts the RKHS embeddings of distributions in the ambiguity set \pfeas.
In this paper, we take a geometric perspective to construct $\mathcal C$ using convex sets in $\mathcal H$.
Given data samples $\{\xi_i\}_{i=1}^N$, we outline various choices for $\mathcal C$ in the left column of Table~\ref{tbl:part} (and \ref{tbl:sets} in the appendix), and illustrate our intuition in Figure~\ref{fig:3in1}.
\begin{figure}[tb!]
        \centering
        \subfloat[RKHS ambiguity sets $\mathcal{C}$\label{fig:3in1_C}]{
        \includegraphics[width=0.83\columnwidth]{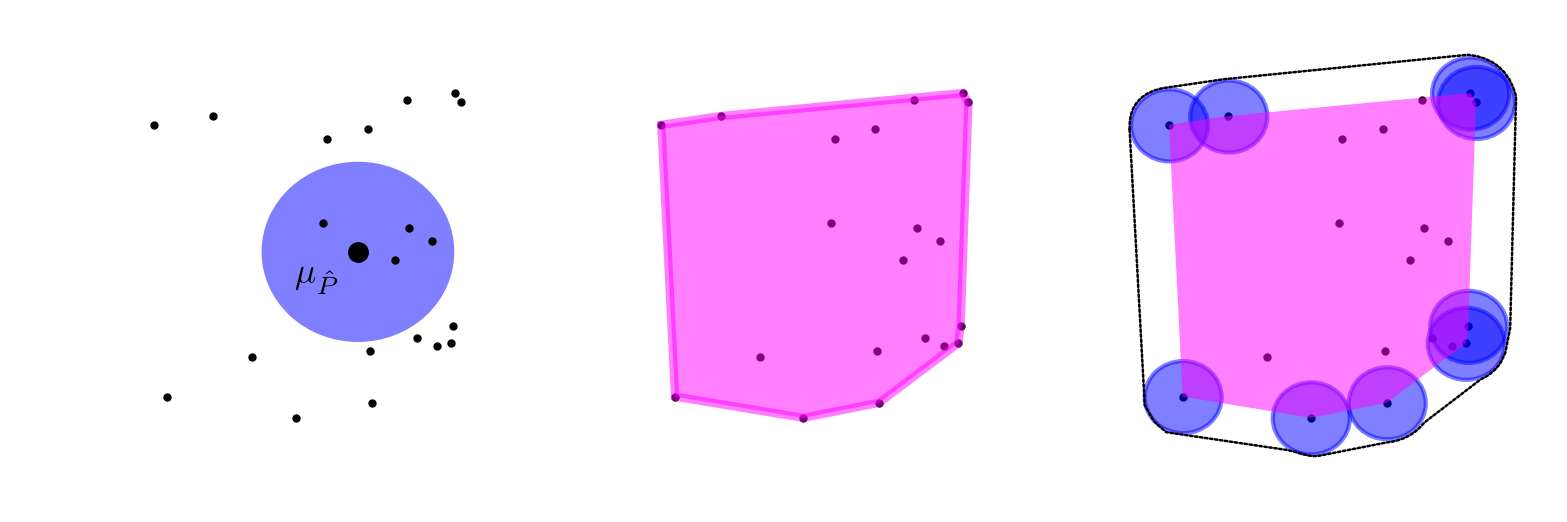}
        }\\
        \subfloat[Interpretation of \kdro\label{fig:3in1_kdro}]{
        \includegraphics[width=0.7\columnwidth]{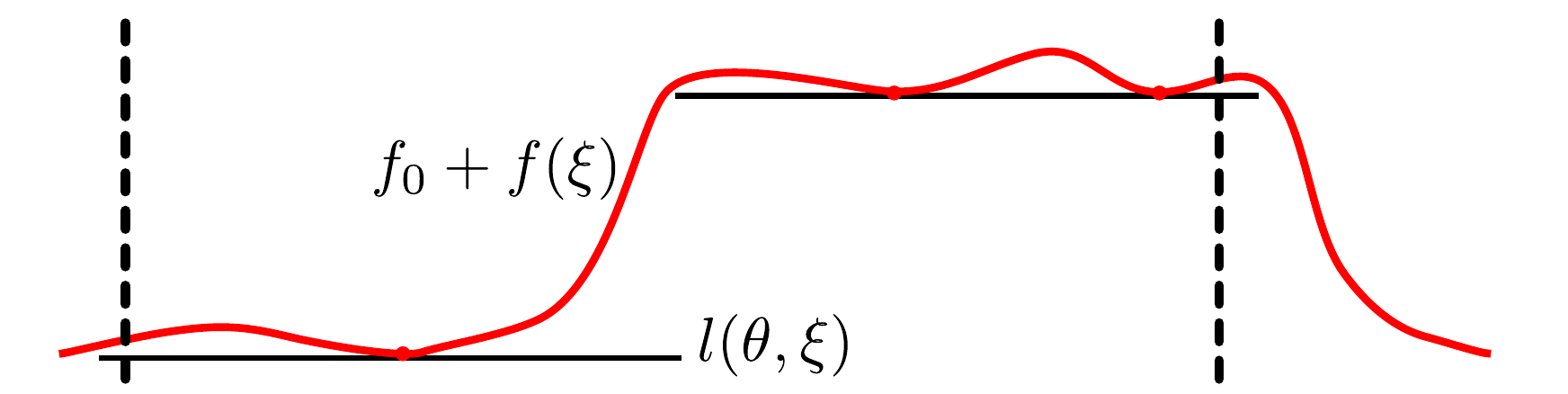}
        }
        \caption{\textbf{(a)}: Geometric intuition for choosing ambiguity set $\mathcal C$ in $\rkhs$ such as
        norm-ball,
        polytope,
        and
        Minkowski sum of sets.
        The scattered points are the embeddings of empirical samples.
        See Table~\ref{tbl:sets} for more examples.
        \textbf{(b)}: Geometric interpretation of \kdro~\eqref{eq:kdro_rc}. The (red) curve depicts $f_0+f$, which \emph{majorizes} $l(\var,\cdot)$ (black). The horizontal axis is $\xi$. The dashed lines denote the boundary of the domain $\mathcal X$.
        }
        \label{fig:3in1}
\end{figure}
\input{draft/dual_various.Rmd}
\input{draft/main_thm.tex}
We now establish further theoretical results as a consequence of the generalized duality theorem to help us understand the geometric intuition of how \kdro works.
By the weak duality $(P)\leq (D)$ of \eqref{eq:moment_primal} and \eqref{eq:dual_conjugate}, we have 
${\int l \ dP  \leq f_0 +\suppc (f)}$.
Specifically, if $\mathcal C $ is the RKHS norm-ball in Table~\ref{tbl:part}
, 
this inequality becomes 
$
        {\int l \ dP \leq f_0+\frac1N\sum_{i=1}^N f(\xi_i)+\epsilon\hnorm{f}}.
$
Its right-hand-side can be seen as a computable bound for the worst-case risk when generalizing to $P$. This may be useful when the Lipschitz constant of $l$ is not known or hard to obtain, as is often the case in practice. 
The following insight is a consequence of a generalization of the classical \emph{complementarity condition} of convex optimization; see the appendix.
\begin{corol}
        [Interpolation property]
        \label{thm:rkhskiss}
        Given \var, let $P^*, f^*, f^*_0$ be a set of optimal primal-dual solutions associated with (P) and (D), then
        $l(\var,\xi) = f_0^* + f^* (\xi)$ holds $P^*$-almost everywhere.
\end{corol}
        Intuitively, this result states that
        $f_0^* +f^*$ interpolates the loss $l(\var,\cdot)$ at the support points of $P^*$.
        This is illustrated in Figure~\ref{fig:3in1} (b) and later empirically validated in Figure~\ref{fig:exp_lsq}.
        We can also see that the size of RKHS $\rkhs$ matters since, if $\rkhs$ is small (e.g., $\rkhs=\{0\}$), $f_0^* +f^*$ cannot interpolate the loss $l$ well. On the other hand, the density of universal RKHS allows the interpolation of general loss functions.

It is tempting to approximately solve \eqref{eq:kdro_rc} by relaxing the constraint to hold for only the empirical samples, i.e., ${l(\var,\xi_i) \leq f_0 + f(\xi_i),
\ i=1\dots N.}$
The following observation cautions us against this.
\begin{ex}[Counterexample: relaxation of the semi-infinite constraint]
        \label{ex:dirac2}
        Let $\rkhs$ be a Gaussian RKHS
        with the bandwidth $\sigma=\sqrt2$. 
        Suppose our data set is $\{0\}$ and the ambiguity set is ${\mathcal C:=\{\mu\colon\|\mu - \phi(0)\|_\rkhs\leq \epsilon \}}$.
        Let $\epsilon=\sqrt{2-2/e}$. 
        We consider the loss function 
        $l(\xi) = [|\var + \xi| - 1]_+$ and relaxing the constraint of \eqref{eq:kdro_rc} to only hold at the empirical sample, i.e.,
\begin{equation*}
(d):=
\begin{opt1}        
        \MINST{\var, f\in\rkhs, f_0\in\mathbb R}{f_0+  f(0) + \epsilon\hnorm{f}}{\sjt\ [|\var | - 1]_+ \leq f_0 + f(0)}
\end{opt1}
\end{equation*}
        which admits an optimal solution $\var^* = 0, f^*=0, f^*_0=0$ and the worst-case risk $(d)=0$. However, let $\mu_{P'}=\frac12\phi( 0) + \frac12\phi( 2)$. It is straightforward to verify $P'\in\mathcal C, \int l(\var^*, \xi) \ dP'(\xi) = \frac12 > (d)$, i.e., the solution $\var^* $ is not robust against $P'$.
\vspace{-0.1cm}
\end{ex}

%% file: draft/dual_various.Rmd
\begin{table*}[t]
	\caption{Examples of support functions for \kdro. See Table~\ref{tbl:sets} for more details. }
	\label{tbl:part}
	\centering
	\begin{tabular}{ll}
	  \toprule
	RKHS ambiguity set $\mathcal C$ & Support function $\delta^*_{\mathcal{C}}(f)$
	\\
	\midrule
	  RKHS norm-ball $\mathcal C=\{\mu\colon\|\mu - \mu_{\hat P}\|_\rkhs\leq \epsilon\}$ &
	  $
	  { \frac1N\sum_{i=1}^N f(\xi_i)+\epsilon\hnorm{f} }
	  $
	  \\
	  Polytope $\mathcal C=\mathrm{conv}\{\phi(\xi_1),\dots,\phi(\xi_N)\}$
	  &
	  ${\max_i f(\xi_i)}$   (scenario optimization; SVMs with no slack)
	  \\
	  Minkowski sum
	  $\mathcal C = \sum_{i=1}^N \mathcal C_i$
	  & 
	  ${ \sum_{i=1}^N\delta^*_{\mathcal C_i} (f)  }$
	  \\
	  Whole space $\mathcal C=\rkhs$ &
	  $0$ if $f=0$, $\infty$ otherwise   (worst-case RO~\cite{ben-talRobustOptimization2009a})
	   \\
	  \bottomrule
	\end{tabular}
\end{table*}

To better understand our unifying formulation, 
let us examine the celebrated SVM through the lens of our generalized formulation.
\begin{ex}[SVM as generalized DRO]
	Let us consider SVM for regression, without using
	slack variables or regularization for simplicity. This can be formulated as optimizing the loss
	$
	\min_{f\in\rkhs}\max_i[|y_i-f(x_i)|-\eta]_+
	$
	where $\eta>0$ is the parameter for the hinge loss. 
	This can be seen as the generalized DRO
	$$
	\min_{f\in\rkhs}\sup_{P\in\mathcal K} \int [|y-f(x)|-\eta]_+ dP(x,y),
	$$
	where the ambiguity set is given by the polytope $\mathcal K = \mathrm{cl}\mathrm{conv}\{\delta_{\xi_1},\dots,\delta_{\xi_N}\}, \xi_i=(x_i,y_i)$. 	
\end{ex}

Let us now consider a small RKHS to understand the effect of different RKHSs on \kdro.
\begin{ex}[DRO with non-universal kernels]
	Consider distributions $\hat{P}=\mathcal{N}(0,1), Q_v=\mathcal{N}(0,v^{2})$ and $\rkhs_1$ induced by the linear kernel $k_{1}(x,y):=xy$. 
	$\rkhs_1$ is small since it only contains linear functions.
	$\mathrm{MMD}_{k_1}(\hat{P},Q_v)=0,\forall v\neq 0$ since they share the first moment.
	Therefore, any $\mathcal C$ that contains $\mu_{\hat P}$ also contains all the distributions in $\{\mu_{Q_v}, v\neq 0 \}$. 
\end{ex}
This example shows that small RKHSs force \kdro to robustify against a large set of distributions, resulting in conservativeness.
In the extreme, if we choose the smallest possible RKHS
$\rkhs=\{0\}$, then the space does not contain functions
to separate any distinct distributions. 
This renders \kdro~\eqref{eq:kdro_rc} overly conservative since we can only choose $f=0$ in \eqref{eq:kdro_rc} ---  it is precisely reduced to worst-case RO.
On the other extreme, the next example shows the downside of function spaces that are too large.
\begin{ex}[DRO with large function space]
Suppose \rkhs is the space of all bounded measurable functions, the metric induced by $\rkhs$ becomes the \emph{total variation} distance~\cite{sriperumbudurUniversalityCharacteristicKernels2011}. While the induced topology is strong, \eqref{eq:kdro_rc} has a trivial solution $f=l,f_0=0$. By plugging this solution into \eqref{eq:kdro_rc}, we recover \eqref{eq:kdro}.	Hence the reformulation becomes meaningless.
\end{ex}

We distinguish between DRO without metrics, e.g., moment constraints and SVMs, and DRO with probability metric or divergence, e.g., Wasserstein metric.
Let us first examine an instance of the former using \kdro. We return to the latter at the end of this section.
\begin{ex}[Reduction to DRO with moment constraints]
	\kdro with the second-order polynomial kernel $k_{2}(x,y):=(1+x^{\top}y)^{2}$ and a singleton ambiguity set $\mathcal C=\{\mu_{\hat{P}}\}$ robustifies 
  against all distributions sharing the first two moments with $\hat{P}$. This is equivalent to 
	DRO with known first two moments, such as in \cite{ delageDistributionallyRobustOptimization2010,scarfMinmaxSolutionInventory1958}.
	More generally, the choice of the $p$th-order polynomial kernel $k_{p}(x,y):=(1+x^{\top}y)^{p}$ corresponds to DRO with known first $p$ moments.
\end{ex}

%
If $\rkhs$ is associated with a universal kernel (e.g., Gaussian), it is large since $\rkhs$ is dense in the space of continuous functions (cf. \cite{sriperumbudurUniversalityCharacteristicKernels2011}). 
Then the induced topology (MMD) is strong enough to separate distinct probability measures. Meanwhile, RKHS allows for efficient computation using tools from kernel methods, as shown in Section~\ref{sec:practical}.
Therefore, our insight is that \emph{universal RKHSs are large enough} for DRO applications.

\begin{remark}
	For the RKHS associated with the Gaussian kernel, the diameter of the space can be computed: $\forall p, q, \|\mu_p - \mu_q\|_\rkhs \leq \|\mu_p \|_\rkhs+\| \mu_q\|_\rkhs
	\leq 2 \sup_{x,y}\sqrt{k(x,y)}= 2$. 
	Hence, if $\epsilon \geq 2$, $\mathcal C$ contains all probability distributions. Then, \kdro is again reduced to worst-case RO on domain \domain.
\end{remark}

We now turn to DRO with a generalized class of integral probability metrics (IPM). 
\begin{ex}[Generalization to IPM-DRO]
	Suppose $\ipm (P,\hat{P}):=\sup_{f\in \mathcal F}\int f d (P-\hat{P})$ is the IPM defined by some function class $\mathcal F$. The IPM-DRO primal formulation is given by

\begin{equation}
\label{eq:ipm_primal}
\min_\var \sup_{\ipm{(P, \hat{P})}\leq\epsilon}\int l(\var, \xi) \ dP(\xi).
\end{equation}
\end{ex}
If we choose the class $\mathcal F=\{f:\hnorm{f}\leq 1\}$, we recover \kdro with the RKHS norm-ball set in Table~\ref{tbl:part}.
Similarly, $\mathcal F=\{f:\mathrm{lip}{(f)}\leq 1\}$ recovers the (type-1) Wasserstein-DRO. 
This puts Wasserstein-DRO and \kdro into a unified perspective.

%% file: draft/main_thm.tex
\subsection{Generalized duality theorem}
\label{sec:dual}
We now present the main theorem of this paper, the generalized duality theorem of \kdro~\eqref{eq:kdro}.

\begin{figure}[h!]
	\vspace{-0.2cm}
	\centering
	\includegraphics[width=0.6\columnwidth]{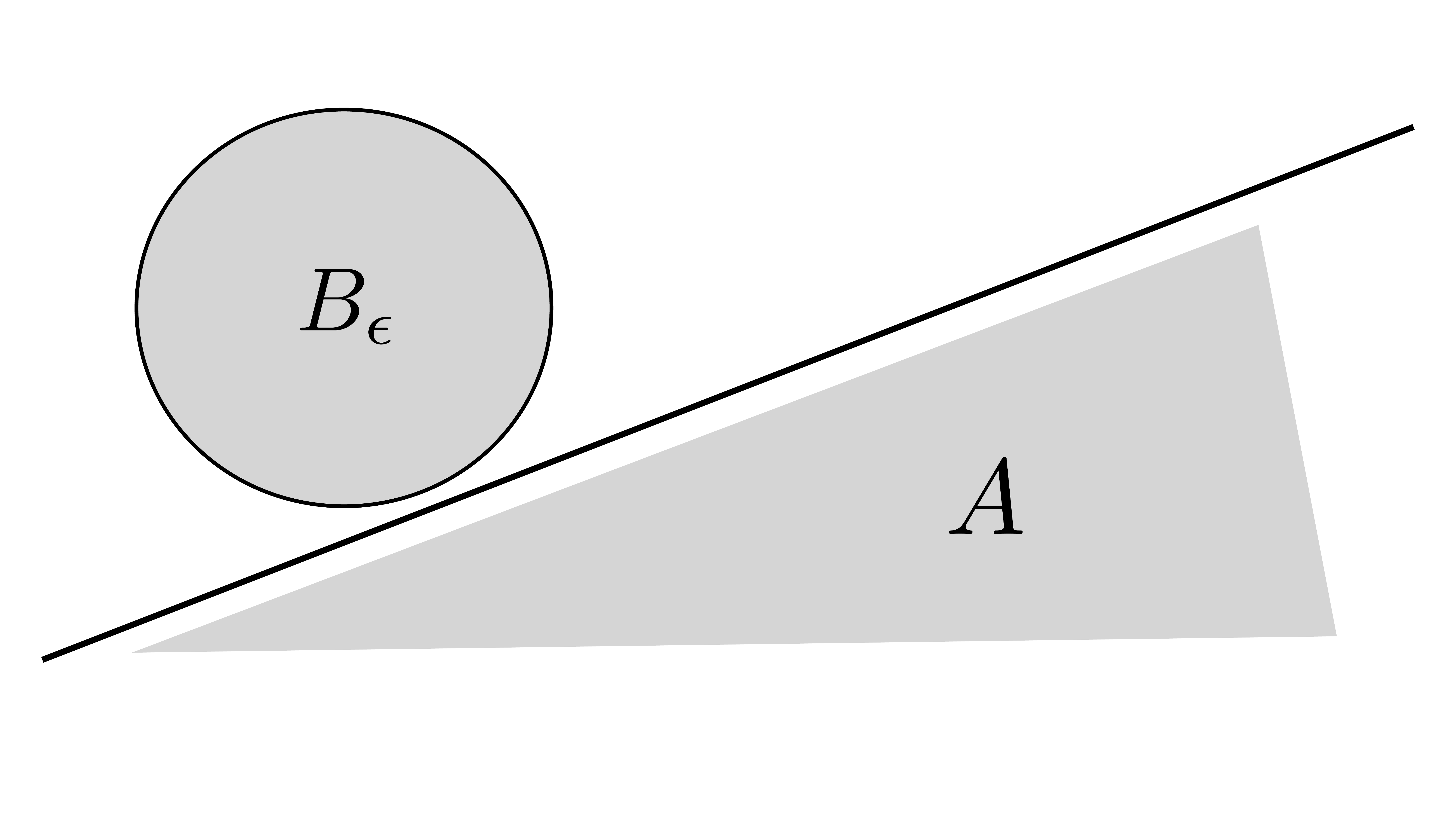}
	\vspace{-0.5cm}
	\caption{Illustration of a separating hyperplane in \rkhs}
	\label{fig:separate_main}
	\vspace{-0.5cm}
  \end{figure}

\begin{theorem}[\textbf{Generalized Duality}]
\label{thm:dro_main_thm}
Under Assumption~\ref{asm:compact}, \eqref{eq:kdro} is equivalent
to
\begin{equation}
	\begin{aligned}
		(D):=&\min_{ \var,f_0\in\mathbb R, f\in\rkhs}&&{\  f_0 + \suppc (f)  }\\
		&\sjt& &\semiinfcons
	\end{aligned}
	\label{eq:kdro_rc}
\end{equation}
where $\suppc (f)  :=\sup_{\mu\in\mathcal C} \hip{f}{\mu} $ is the support function of $\mathcal C$, i.e., $(P)=(D)$, \emph{strong duality} holds for the inner moment problem for any \var point-wise.
\end{theorem}
The theorem holds regardless of the dependency of $l$ on \var, e.g., non-convexity.
If $l$ is convex in \var, then \eqref{eq:kdro_rc} is a \emph{convex program}.
Formulation \eqref{eq:kdro_rc} has a clear geometric interpretation: we find a function $f_0+f$ that \emph{majorizes} $l(\var,\cdot)$ and subsequently minimize a surrogate loss involving $f_0$ and $f$. This is illustrated in Figure~\ref{fig:3in1_kdro}. Note the term duality here refers to the inner moment problem.
The statement can be further simplified by replacing $f_0+f$ with $f$. However, we choose the current notation for the sake of its explicit connection to RO.
\paragraph*{Proof sketch.}
Our weak duality proof follows standard paradigms of Lagrangian relaxation by introduing dual variables.
Notably, we associate the functional constraint $\int{\phi }{\ dP} = \mu$ with a dual function $f\in \rkhs$, which is the decision variable in the dual problem~\eqref{eq:kdro_rc}.
Using the reproducing property of RKHSs and conic duality, we arrive at~\eqref{eq:kdro_rc} with weak duality.
Our strong duality proof is an extension of the conic strong duality in Eulidean spaces.
We rely on the existance of separating hyperplnes between convex sets in locally convex function spaces, e.g., \rkhs. See the illustration in Figure~\ref{fig:separate_main}. In our generalized duality theorem, this separating hyperplane is determined by the witness function $f^*$, which is the optimal dual variable in \eqref{eq:kdro_rc}. See the appendix for the full proof.

Theorem~\ref{thm:dro_main_thm} generalizes the classical bounds in generalized moment problems \cite{isiiSharpnessTchebychefftypeInequalities1962,lasserreBoundsMeasuresSatisfying2002,
shapiroDualityTheoryConic2001,
bertsimasOptimalInequalitiesProbability2005,popescuSemidefiniteProgrammingApproach2005,vandenbergheGeneralizedChebyshevBounds2007,
vanparysGeneralizedGaussInequalities2016} to infinitely many moments using RKHSs.
A distinction between Theorem~\ref{thm:dro_main_thm} and other DRO approaches is that it uses the density of universal RKHSs to find a surrogate which can sharply bound the worst-case risk. This means that we do not require the loss $l(\var,\cdot)$ to be affine, quadratic, or living in a known RKHS, nor do we require the knowledge of Lipschitz constant or RKHS norm of $l(\var,\cdot)$.
To our knowledge, existing works
typically require one of such assumptions.

Moreover, Theorem~\ref{thm:dro_main_thm} generalizes existing RO and DRO in the sense that it gives us a unifying tool to work with various ambiguity and ambiguity sets, which may be customized for specific applications.
We outline a few closed-form expressions of the support function $\suppc (f)$ in Table~\ref{tbl:part}, and more in Table~\ref{tbl:sets}. 
We now return IPM-DRO with a duality result.
\begin{corol}[IPM-DRO duality]
	\label{thm:ipm_dual}
Given the integral probability metric $\ipm (P,\hat{P}):=\sup_{f\in \mathcal F}\int f d (P-\hat{P})$, a dual program to \eqref{eq:ipm_primal} is given by
	\begin{equation}
		\label{eq:ipm0}    
		\begin{opt1}
			\MINST{\var, \lambda\geq 0, f_0\in\mathbb R, f\in \mathcal F}{\           {f_0+\frac{1}N\sum_{i=1}^N \lambda f(\xi_i)+\lambda\epsilon }
			}{l(\var, \xi)\leq  f_0 + \lambda f(\xi)   ,\ \forall \xi \in\mathcal X.} 
		\end{opt1}
	\end{equation}
\end{corol}
The reduction to \eqref{eq:kdro_rc} as a special case can be seen by replacing $\lambda f$ with $f$ and choosing $\mathcal F=\rkhs$. 


We now establish an explicit connection between DRO and stochastic optimization with expectation constraint, whose solution methods using stochastic approximation are an topic of active research~\cite{lanAlgorithmsStochasticOptimization2020a,xuPrimalDualStochasticGradient2020}.
\begin{corol}(Stochastic optimization with expectation constraint)
	\label{thm:expect_cons}
	Under the Assumption~\ref{asm:compact}, the optimal value of \eqref{eq:kdro} coincides with that of
	\begin{equation}
		\begin{opt1}
		\label{eq:sa_exp_cons}
		\MINST{ \var,f_0\in\mathbb R, f\in\rkhs}{\  f_0 + \suppc (f)  }
		{\mathbb{E}h\left(l\left(\var, \zeta\right) -   f_0 - \lambda f\left(\zeta\right) \right) \leq 0}
		\end{opt1}
	  \end{equation}
	  for some function $h$ that satisfies $h(t) = 0 \text{ if } t\leq 0, h(t)>0 \text{ if } t>0$, and random variable $\zeta\sim\mu$ whose probability measure places positive mass on any non-empty open subset of \domain, i.e.,
	  $
	  \mu (B)>0,\ \forall B\subseteq \domain, B\neq\emptyset, B\text{ is open}
	  $.
\end{corol}
A choice for $h$ is $h(\cdot) = [\cdot]_+$, which is used in the conditional value-at-risk~\cite{rockafellarOptimizationConditionalValueatrisk2000a}.
We will see the computational implication of Corollary~\ref{thm:expect_cons}
in Section~\ref{sec:practical}.

%% file: draft/practical.Rmd

Given a certain parametrization of the RKHS function $f$, \eqref{eq:kdro_rc} is a semi-infinite program  (SIP)~\cite{guerravazquezGeneralizedSemiinfiniteProgramming2008}. 
In the following, we propose two computational methods
that do not require a polynomial loss $l$ or the knowledge of its Lipschitz constant.
For simplicity, we only derive the formulations for the RKHS-norm-ball ambiguity set, while other formulations are given in Table~\ref{tbl:part},~\ref{tbl:sets}.


\paragraph*{A batch approach by discretization of SIP.}
We first consider an approach based on the discretization method of SIP~\cite{guerravazquezGeneralizedSemiinfiniteProgramming2008}.
Let us consider an ambiguity set smaller than the RKHS-norm ball of distributions supported on some $\{\zeta_j\}_{j=1}^M\subseteq\mathcal X$.
Then it suffices to consider the following program, which relaxes the constraint of \eqref{eq:kdro_rc} to finite support.
\begin{equation}
                \begin{opt1}
                \MINST{\var, f\in\rkhs, f_0\in\mathbb R}
                {f_0+\frac1N\sum_{i=1}^N f(\xi_i) + \epsilon\hnorm{f}}
                {l(\var,\zeta_i) \leq f(\zeta_j)  + f_0 ,\  j=1\dots M.}
                \end{opt1}
        \label{eq:scenario}
\end{equation}
Note \eqref{eq:scenario} is a \emph{convex program} if $l(\var,\xi)$ in convex in $\var$.

We can parametrize the RKHS function $f$ by a wealth of tools from kernel methods, such as the random features $\hat{f}(\xi)= w^\top \hat\phi(\xi)$ for large scale problems~\cite{rahimiRandomFeaturesLargeScale2008}.
Alternatively, for small problems, we can parametrize $f$ by a kernel expansion on the the points $\zeta_i$.
We provide concrete plug-in forms in the appendix.

As an interesting by-product of \eqref{eq:scenario}, let us derive  an unconstrained version of \eqref{eq:scenario}, which gives rise to a generalized risk measure that we term \emph{kernel conditional value-at-risk} (Kernel CVaR).
\begin{ex}[Kernel CVaR]
        \label{ex:kcvar}
        \begin{multline}
                \vspace{-2cm}
                        \operatorname{K-CVaR}_\alpha(X):=\inf_{ f\in\rkhs, f_0\in\mathbb R}
                        \ \frac1{\alpha M}\sum_{j=1}^M 
                        [X - f(\zeta_j) - f_0]_+\\
                        + f_0+\frac1N\sum_{i=1}^N f(\xi_i) + \epsilon\hnorm{f}.
                        \label{eq:ridge_cvar}
                \end{multline}
                If $f=0$ and $\{\zeta_j\}_{j=1}^M=\{\xi_i\}_{i=1}^N$, then \eqref{eq:ridge_cvar} is reduced to the classical $\operatorname{CVaR}$~\cite{rockafellarOptimizationConditionalValueatrisk2000a}.
        \end{ex}
Program~\eqref{eq:scenario} can be readily solved using off-the-shelf convex solvers.
However, to scale up to large data sets, we next develop a stochastic approximation (SA) method for \kdro.

%% file: draft/aistats/sgd_kdro.tex
\paragraph*{Stochastic functional gradient DRO.}
We now present our SA approach enabled by Theorem~\ref{thm:dro_main_thm}
by employing two key tools:
1) scalable approximate RKHS features, such as random Fourier features~\cite{rahimiRandomFeaturesLargeScale2008,daiScalableKernelMethods2014,carratinoLearningSGDRandom2018}, and
2) stochastic approximation with semi-infinite and expectation constraints~\cite{tadicRandomizedAlgorithmsSemiInfinite2006,lanAlgorithmsStochasticOptimization2020a,baesRandomizedMirrorProxMethod2011,xuPrimalDualStochasticGradient2020}.

Let us summon Corollary~\ref{thm:expect_cons} to formulate a stochastic program with expectation constraint.
\begin{equation}
  \begin{opt1}
  \label{eq:sa}
  \MINST{ \var,f_0\in\mathbb R, f\in\rkhs}{\  f_0+\frac1N\sum_{i=1}^N f(\xi_i) + \epsilon\hnorm{f}  }{\mathbb{E} [l\left(\var, \zeta\right) -   f_0 - \lambda f\left(\zeta\right) ]_+ \leq 0}
  \end{opt1}
\end{equation}
where $\zeta$ follows a certain proposal distribution on \domain, e.g., uniform or by adaptive sampling. 
An alternative is to directly solve \eqref{eq:kdro_rc} using SA techniques with SI constraints, such as \cite{tadicRandomizedAlgorithmsSemiInfinite2006,weiCoMirrorAlgorithmRandom2020a}.
\eqref{eq:kdro_rc}, \eqref{eq:sa_exp_cons}, and \eqref{eq:sa} are all convex in function $f$.
We can compute the functional gradient
by
\begin{equation}
\label{eq:fgrad}
\nabla_f f
= \phi, \quad \nabla_f \hnorm{f} =\frac{f}{\|f\|_\rkhs}.
\end{equation}

When used with approximate features of the form 
$\hat{f}(\xi)= w^\top \hat\phi(\xi)$, 
we further have
$
\nabla_w \hat{f}(\xi) = \hat\phi(\xi), \nabla_w \hnorm{\hat f} = w/\|w\|_2.
$
We outline our stochastic functional gradient DRO (\sfg) in Algorithm~\ref{alg:sgd}.
\begin{algorithm}[htb!]
  \caption{Stochastic Functional Gradient DRO (\sfg)}
  \label{alg:sgd}
  \begin{algorithmic}[1]
      \FOR {$k=1,2,\dots$}
      \STATE Sample mini-batch data $\{\xi^k_i\}:=\{x_i, y_i\}$. Sample $\{\zeta_i\}$ from some proposing distribution. 
      \label{step:sample}
      \STATE Approximate $f$, e.g., by random Fourier feature $\hat{f}(\xi^k_i)= w^\top \hat\phi(\xi^k_i)$.
      \label{step:feature}
      \STATE Estimate the stochastic functional gradient of the objective and constraint in \eqref{eq:sa} using \eqref{eq:fgrad}.
      \label{step:grad}
      \STATE Update $\var, f_0, f$ using the functional gradient with any SA routine with expectation or semi-infinite constraints, e.g., \cite{lanAlgorithmsStochasticOptimization2020a,xuPrimalDualStochasticGradient2020,tadicRandomizedAlgorithmsSemiInfinite2006,weiCoMirrorAlgorithmRandom2020a}
      \label{step:sa}.
      \ENDFOR
  \end{algorithmic}
  \end{algorithm}

Compared with many batch-setting DRO approaches, \sfg can be used with general model classes, such as neural networks, and is applicable to a broad class of optimization and modern learning tasks. The convergence guarantee follows that of the specific SA routine used in Step~\ref{step:sa} of the algorithm. It is worth noting that, when used with a primal SA approach such as \cite{lanAlgorithmsStochasticOptimization2020a}, \sfg completely operates in the dual space (an RKHS) since \kdro~\eqref{eq:kdro_rc} is based on the generalized duality Theorem~\ref{thm:dro_main_thm}. This interplay between the primal (measures) and dual (functions) is the essence of our theory.

%% file: draft/aistats/exp_fig_all.tex
\begin{figure*}[h]
	\centering
	\subfloat[Uncertain least squares loss\label{fig:lsq_test}]{
	\includegraphics[height=3.6cm,valign=c]{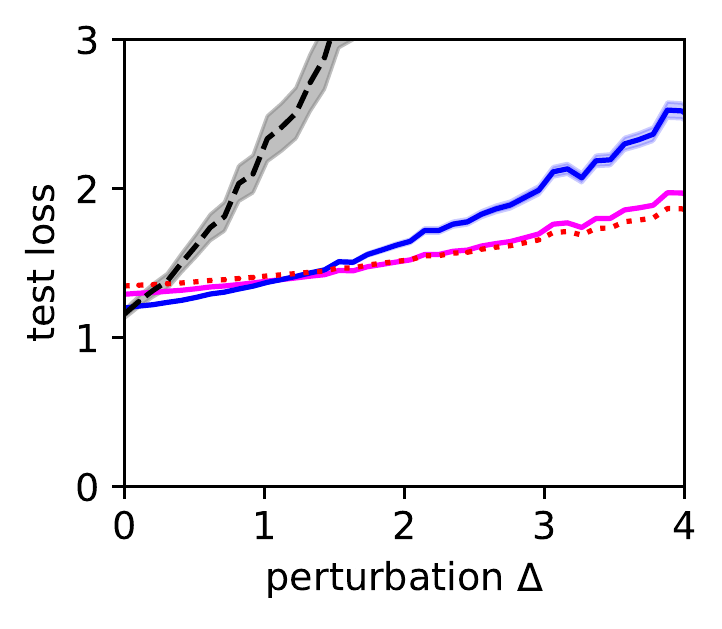}
	\includegraphics[width=1.5cm,valign=c]{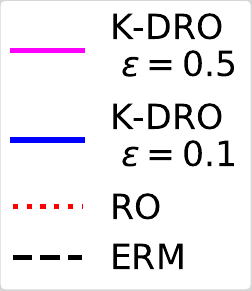}
	}
	\subfloat[Geometric interpretation\label{fig:lsq_interpret}]{
	\includegraphics[height=3.7cm,valign=c]{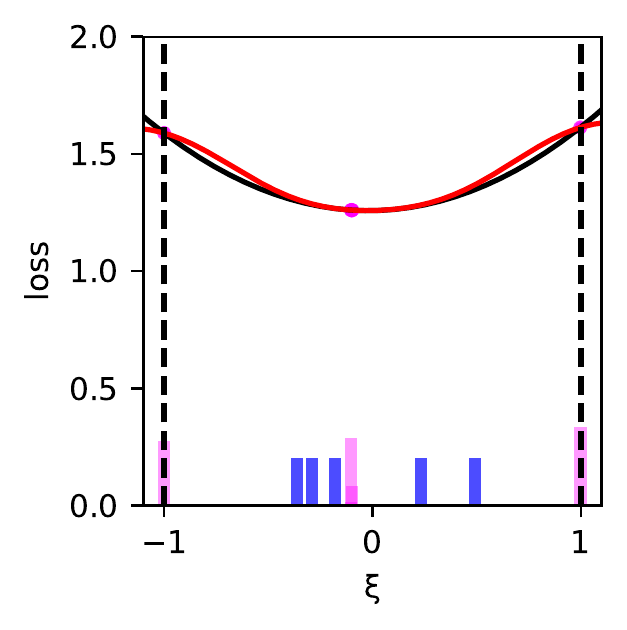}
	\includegraphics[width=1.4cm,valign=c]{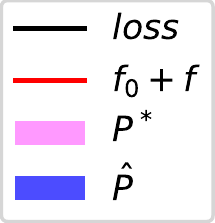}
	}
	\subfloat[MNIST classification error\label{fig:mnist_curve}]{
        \includegraphics[height=3.7cm,valign=c]{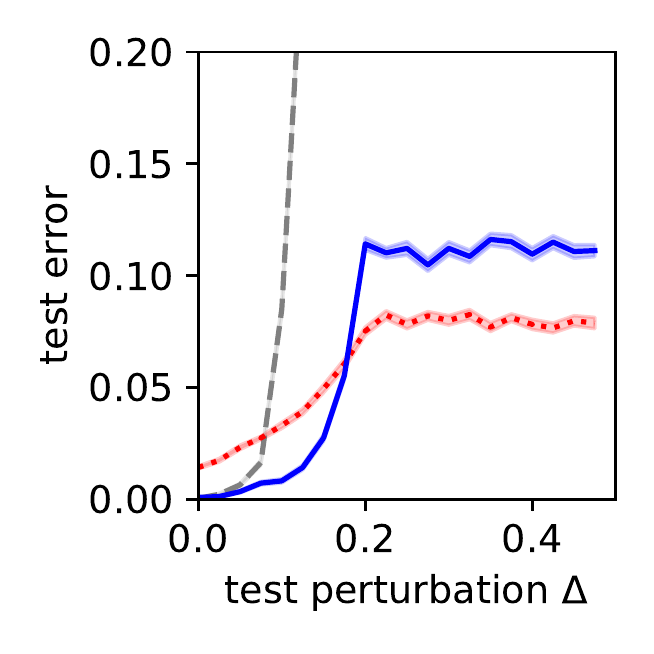}
        \includegraphics[width=1.7cm,valign=c]{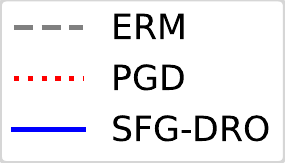}
        }
	\par
	\subfloat[(Left) unperturbed data (Center) ERM classification result (red indicates errors) (Right) \sfg (\kdro)\label{fig:mnist_viz}]{
            \includegraphics[height=1.6cm,valign=c]{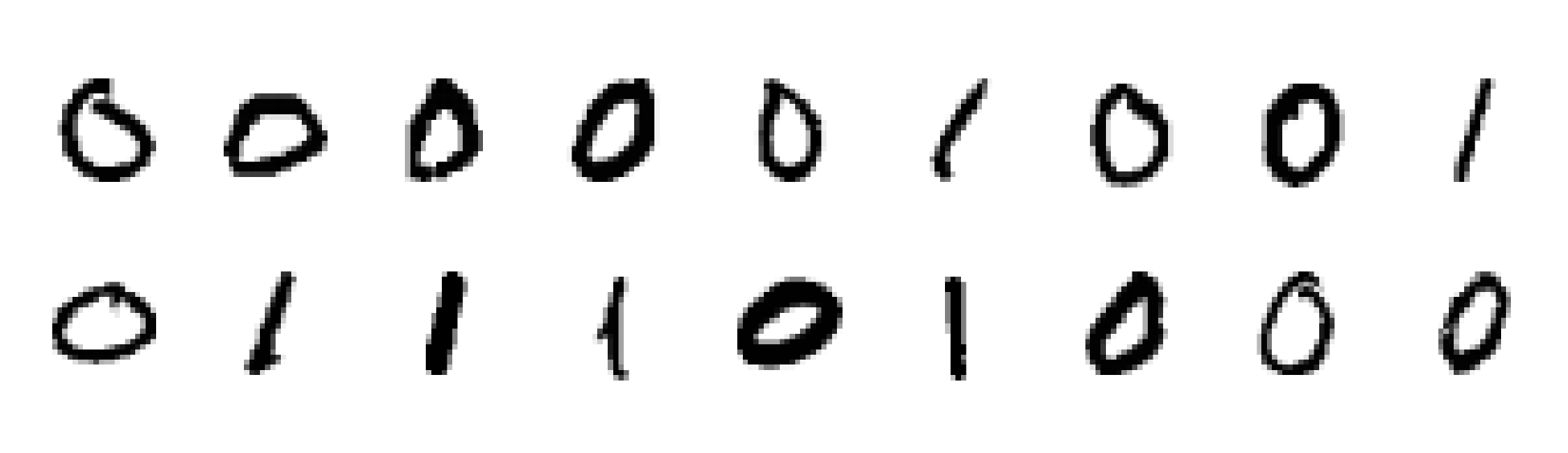}\par
            \includegraphics[height=1.6cm,valign=c]{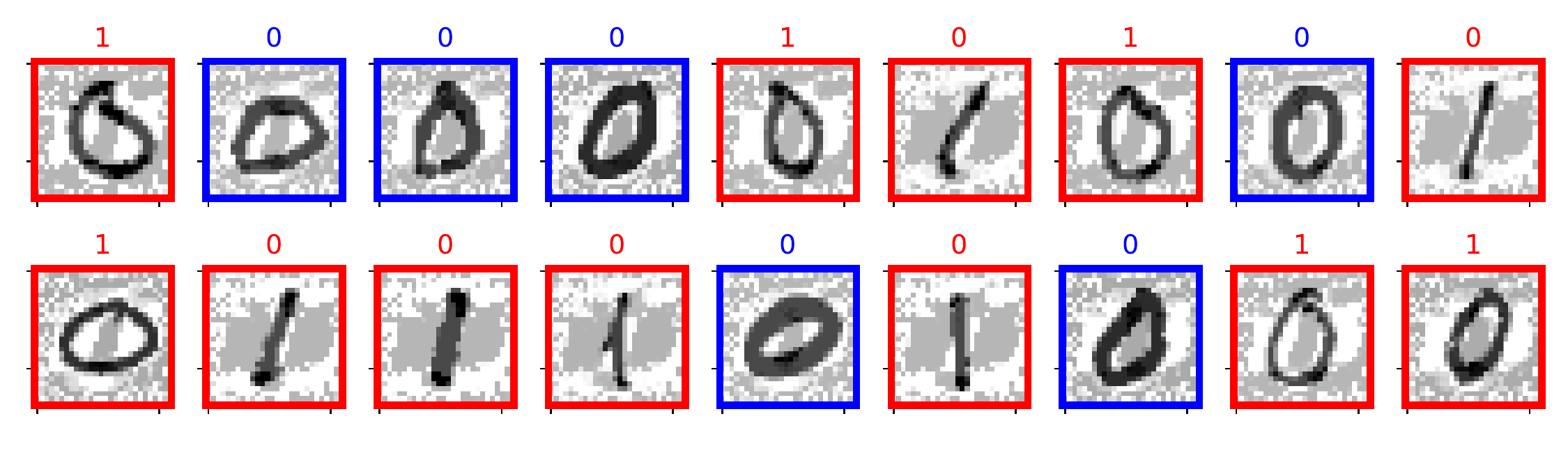}\par
            \includegraphics[height=1.6cm,valign=c]{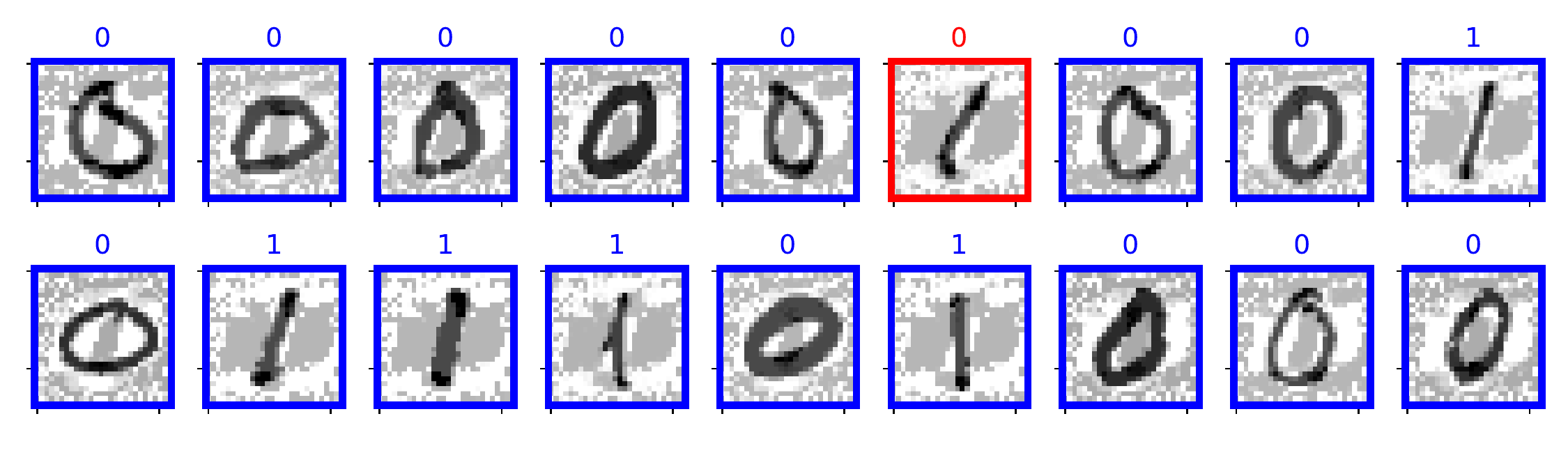}
    }
	\caption{
		\textbf{Uncertain least squares.}
		\textbf{(a)} This plot depicts the test loss of algorithms. 
		All error bars are in standard error.
		We ran $10$ independent trials. In each trial, we solved \kdro to obtain \var$^*$ and tested it on a test dataset of 500 samples. We then vary the perturbation $\Delta$ from $0$ to $4$.
		\textbf{(b)}
		(red) is the dual optimal solution $f_0^*+f^*$.
		(black) is the function $l(\var ^*, \cdot)$.
		The pink bars depict a worst-case distribution while the blue bars the empirical distribution. 
		We can observe that $f_0^*+f^*$ touches loss $l(\var ^*, \cdot)$ at the support of the worst-case distribution $P^*$ (pink dots).
		Note ${f^*}$ (normalized) can be viewed as a witness function of the two distributions.
		\textbf{Classification under perturbation}
        \textbf{(c)}
		We plot the classification error rate during test time.
		The $x-$axis is the perturbation magnitude allowed on the test data.
		For ERM, PGD, and \sfg (\kdro), we train $5$ independent models. Each model is tested on $500$ randomly sampled images.
		\textbf{(d)}
		We visualize the predictions of ERM and \sfg on the perturbed images with perturbation magnitude $\Delta=0.2$. Blue frames indicate correct predictions while the red ones indicate errors.
		}
		\label{fig:exp_lsq}
		\vspace{-2mm}
\end{figure*}

%% file: draft/exp_1.tex
We first consider a robust least squares problem adapted from \cite{elghaouiRobustSolutionsLeastSquares1997}, which demonstrated an important application of RO to statistical learning historically. (See also \cite[Ch.~6.4]{boydConvexOptimization2004}.)
The task is to minimize the objective $\|A\var-b\|_2^2$ w.r.t. $\var$. 
$A$
is modeled by $A(\xi) = A_0 + \xi A_1$, where 
$\xi\in\domain$ is uncertain,
$\domain=[-1,1]$,
and
$A_0, A_1 \in\mathbb R^{10\times 10}, b\in\mathbb R^{10}$ are given.
We compare \kdro against using \emph{(a)} empirical risk minimization (ERM; also known as sample average approximation) that minimizes $\frac1N\sum_{i=1}^N {\|A(\xi_i)\ \var-b\|_2^2}$, \emph{(b)} worst-case RO via SDP from \cite{elghaouiRobustSolutionsLeastSquares1997}.
We consider a data-driven setting with given samples $\{\xi_i\}_{i=1}^N$
with the \kdro formulation
$
\mnz _ \var \max_{P\in\mathcal P,\mu\in\mathcal C}  \expected{\xi\sim P}{\|A(\xi)\ \var-b\|_2^2}\quad
\sjt \int \phi dP = \mu,
$
where we choose the ambiguity set to be
the $\epsilon$-norm-ball in the RKHS (Table~\ref{tbl:part}).



Empirical samples $\{\xi_i\}_{i=1}^N (N=10)$ are generated uniformly from $[-0.5, 0.5]$. We then apply \kdro formulation~\eqref{eq:scenario}.
To test the solution, we create a distribution shift by generating test samples from $[-0.5\cdot(1+\Delta), 0.5\cdot(1+\Delta)]$, where $\Delta$ is
a perturbation varying within $[0,4]$.
Figure~\ref{fig:lsq_test} shows this comparison.
As the perturbation increases, ERM quickly lost robustness. 
On the other hand, RO is the most robust 
with the trade-off of being conservative. 
As expected, \kdro achieves some level of optimality while retaining robustness. 

We then ran \kdro with fewer empirical samples ($N=5$) to show the geometric interpretations.
We plot the optimal dual solution $f^*_0 + f^*$ in Figure~\ref{fig:lsq_interpret}. Recall it is an over-estimator of the loss $l(\var, \cdot)$.
We solve the inner moment problem (see appendix) to obtain a worst-case distribution $P^*$. 
Comparing $P^*$ with $\hat P$, 
we can observe the adversarial behavior of the worst-case distribution.
See the caption for more description.
From Figure~\ref{fig:lsq_interpret}, we can see that the \emph{intuition of \kdro is to flatten the loss curve using a smooth function}.

%% file: draft/aistats/exp_adv.tex
\paragraph{Distributionally robust learning under adversarial perturbation.}
We now demonstrate the framework of \sfg in Algorithm~\ref{alg:sgd} in a non-convex setting.
For simplicity, we consider a MNIST binary classification task with a two-layer neural network. 
We emphasize that the deliberate choice of this simple architecture ablates factors known to implicitly influence robustness, such as regularization and dropout.
The data set contains images of zero and one (i.e., two classes).
Each image $x$ is represented by $x\in[0,1]^{28 \times 28}$.
The test data is perturbed by an \emph{unknown} disturbance, i.e., $\tilde x_\mathrm{test} : = x + \delta$ where $x\sim P_{\mathrm{test}}$ is the unperturbed test data and $\delta$ is the perturbation.
In the plots, $\delta$ is generated by the PGD algorithm~\cite{madryDeepLearningModels2019} using projected gradient descent to find the worst-case perturbation within a box $\{\delta:\|\delta\|_\infty\leq \Delta\}\}$.
We compared \sfg (\kdro) with ERM and PGD (cf. \cite{madryDeepLearningModels2019,madryAdversarialRobustnessTheory}).
Note
the overall loss of PGD is an average loss instead of a worst-case one. Hence it is already less conservative than RO. 
We train a classification model $g_\var\colon x \mapsto y$ using SFG-DRO in Algorithm~\ref{alg:sgd}, with the SA subroutine of \cite{lanAlgorithmsStochasticOptimization2020a}.
During training, we set the ambiguity size of \sfg as $\epsilon=0.5$ and domain \domain to be norm-balls around the training data $\domain=\{\zeta = X + \delta: \|\delta\|_\infty\leq 0.5\}$ where $X$ is the training data. 

Figure~\ref{fig:mnist_viz}~(left) plots unperturbed test samples.
Figure~\ref{fig:mnist_curve} shows the classification error rate as we increase the magnitude of the perturbation $\Delta$.
We observe that ERM attains good performance when there is no test-time perturbation but quickly underperforms as the noise level increases.
PGD
is the most robust under large perturbation, but has the worst nominal performance.
\sfg possesses improved robustness while its performance under no perturbation does not become much worse. 
This is consistent with our theoretical insights into RO and DRO.

%% file: draft/literature.tex
This paper uses similar techniques of reformulating min-max programs as in \cite{ben-talDerivingRobustCounterparts2015,bertsimasDatadrivenRobustOptimization2017}, but our ambiguity set is constructed in an RKHS.
\citet{duchiDistributionallyRobustLosses}
proposed variational approximations to marginal DRO to treat covariate shift in supervised learning.
The authors of \cite{zhuWorstCaseRiskQuantification2020} used kernel mean embedding for the inner moment problem (but not DRO) and proved the statistical consistency of the solution.
The work of \citet{StaJeg2019} used insights from DRO to motivate a regularizer for kernel ridge regression.
DRO has been also applied to Bayesian optimization in \cite{rontsisDistributionallyAmbiguousOptimization2020,kirschnerDistributionallyRobustBayesian2020}, where the latter work used MMD ambiguity sets of distributions over discrete spaces.
In terms of scalability, recent works such as \cite{sinhaCertifyingDistributionalRobustness2020a,liFirstOrderAlgorithmicFramework2019,namkoongStochasticGradientMethods2016} also explored DRO for modern machine learning tasks.
To the best of our knowledge, no existing work contains the results
 such as generalized ambiguity set constructions in Table~\ref{tbl:part},~\ref{tbl:sets}, generalized duality theory underpinned by Theorem~\ref{thm:dro_main_thm}, or the stochastic functional gradient algorithm \sfg.

%% file: draft/discuss.Rmd
\emph{In summary}, this paper proves Theorem~\ref{thm:dro_main_thm} that generalizes the classical duality theory in the literature of mathematical problem of moments and DRO. Using the density of universal RKHSs, the dual bound in Theorem~\ref{thm:dro_main_thm} is sharp while lifting restrictions on the loss function class.
The generalized primal formulations shed light on the connection between \kdro and existing robust and stochastic optimization approaches.
Finally, the proposed stochastic approximation algorithm \sfg enables the applications of \kdro to modern learning tasks.

The compactness assumption on \domain can be further extended,
as universality can be extended to non-compact domains~\cite{sriperumbudurUniversalityCharacteristicKernels2011}.
In the special case of RKHS-norm-ball ambiguity sets, choosing the size $\epsilon$ can be motivated using kernel statistical testing~\cite{grettonKernelTwosampleTest2012}. However, when DRO is used in the setting where test distributions are perturbed as in our examples, existing statistical guarantees in the literature for unperturbed settings cannot be directly applied. This is a topic of future work.

%% file: draft/proof.Rmd
Table~\ref{tbl:landscape} provides an overview to help readers navigate the theoretical results in this paper. 
\begin{table}[h]
        \caption{List of the theoretical results in this paper}
        \label{tbl:landscape}
        \centering
        \begin{tabular}{ll}
          \toprule
          Theorems 
          & Generalized Duality Theorem~\ref{thm:dro_main_thm}\\
          & Strong duality of the inner moment problem Proposition~\ref{thm:dual_conjugate}\\
          & Interpolation property Proposition~\ref{thm:rkhskiss}
          \\
          & Complementarity condition Lemma~\ref{thm:complement}\\
          & Robust representer theorem Proposition~\ref{thm:representer_weak}\\ 
          & IPM-DRO duality Corollary~\ref{thm:ipm_dual}\\
          & \kdro as stochastic optimization with expectation constraint Corollary~\ref{thm:expect_cons}\\

          \midrule
          Formulations & \kdro primal (P)~\eqref{eq:kdro}, dual (D)~\eqref{eq:kdro_rc}\\
          & IPM-DRO primal~\eqref{eq:ipm_primal}, dual~\eqref{eq:ipm0}\\
          & Formulations for various RKHS ambiguity sets Table~\ref{tbl:part},~\ref{tbl:sets}\\
          & Stochastic program with expectation constraint formulation of \kdro~\eqref{eq:sa_exp_cons},\eqref{eq:sa}\\
          & Program to compute worst-case distributions~\eqref{eq:cdc},\eqref{eq:perturb},\eqref{eq:sdp}\\
          & \kdro convex program by the discretization of SIP~\eqref{eq:scenario}\\
          & Kernel conditional value-at-risk~\eqref{eq:ridge_cvar}\\
          \bottomrule
        \end{tabular}
\end{table}

In general, we refer to standard texts in optimization \cite{boydConvexOptimization2004,shapiroLecturesStochasticProgramming2014,ben-talRobustOptimization2009a}, convex analysis \cite{rockafellarConvexAnalysis1970,barvinokCourseConvexity2002}, and functional analysis \cite{conway2019course}
for more mathematical background.
\paragraph{Notation.}
In the proofs, we use
$\mathcal M$ to denote the space of signed measures on \domain.
The dual cone of a set of signed measures $\mathcal K\subseteq\mathcal M$ is defined as
$
\mathcal K^*:= \{h\colon\int h \, \mathrm{d}m\geq0, \forall m \in \mathcal K, h\ \ \mathrm{measurable}\}.
$
Using the reproducing property, we have the identity $\pip fP = \hip{f}{\mu_P}$ for $f\in\rkhs$ and $P\in\mathcal P$, which we will frequently use in the proofs.

\subsection{Proof of the Generalized Duality Theorem~\ref{thm:dro_main_thm}}
\begin{figure}[htb!]
        \centering
        \includegraphics[width=0.4\linewidth]{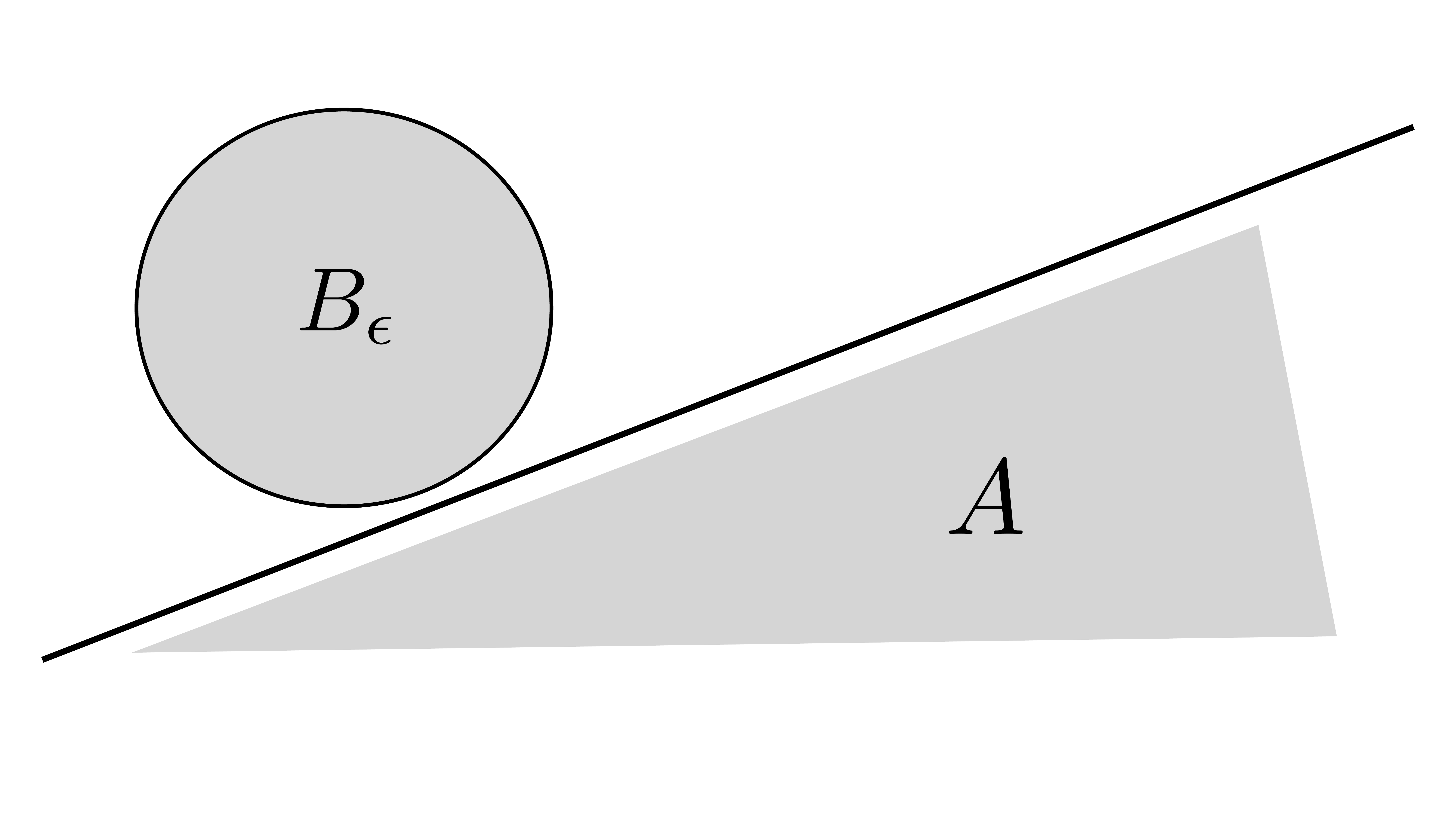}
        \caption{Illustration of the strong duality proof that uses a separating hyperplane. See the proof for detailed descriptions.}
        \label{fig:separate}
      \end{figure}

We now derive our key result for \kdro\ --- the Generalized Duality Theorem, in Theorem~\ref{thm:dro_main_thm}. 
Let us first consider the inner moment problem of \eqref{eq:kdro}
\begin{equation}
        \begin{opt1}
\SUPST{P\in \mathcal P,\mu\in \mathcal C}{\int l \ dP}
{\int{\phi}{\ dP} = \mu},
        \end{opt1}        
\label{eq:moment_primal}\end{equation}
where we suppress \var in $l(\var, \cdot)$ as we fix it for the moment.
\eqref{eq:moment_primal} generalizes the \emph{problem of moments} in the sense that the constraint can be viewed as infinite-order moment constraints.
Using conic duality, we obtain the strong duality of the inner moment problem.
\begin{prop}[Strong dual to \eqref{eq:moment_primal}]
        \label{thm:dual_conjugate}
Under Assumption~\ref{asm:compact},
\eqref{eq:moment_primal} is equivalent to solving
\begin{equation}        
\begin{opt1}
        \MINST{ f_0\in\mathbb R, f\in\rkhs}{\ \suppc  (f) + f_0}
        {l(\xi)\leq  f_0 + f(\xi)   ,\ \forall \xi \in\mathcal X} 
\end{opt1}
\label{eq:dual_conjugate}
\end{equation}
where $\suppc$ is the support function of $\mathcal C$, i.e., \emph{strong duality} holds.
\end{prop}
Using Proposition~\ref{thm:dual_conjugate}, we can reformulate the inner moment problem in \eqref{eq:kdro} to obtain Theorem~\ref{thm:dro_main_thm}.
We now prove this generalized duality result for the inner moment problem in Proposition~\ref{thm:dual_conjugate}.
We first derive the weak dual and then prove the strong duality.
\begin{proof}
We first relax the constraint $P\in \mathcal P$ to its conic hull $P\in \co{\mathcal P}$. To constrain $P$ to still be a probability measure, we impose $\int{1}{\ dP(x)} = 1$, which results in the primal problem equivalent to \eqref{eq:moment_primal}
$$
\begin{opt1}
        (P):=\MAXST{P\in \co{\mathcal K},\mu\in\mathcal C}{\pip lP}{{\int \phi\ dP} = \mu,\ \pip 1P=1}.
\end{opt1}
$$
We construct the Lagrangian relaxation by associating the constraints with the dual variables ${f\in\mathcal H, f_0\in\mathbb R}$, as well as adding the indicator function of $\mathcal C$.
Note both sides of the constraint $\int{\phi }{\ dP} = \mu$ are functions in $\rkhs$, hence the multiplier $f$ is an RKHS function. 
\begin{multline}
\label{eq:lagrangian}
\mathcal L(P, \mu;\  f,f_0) =
\pip lP - \indc (\mu) + \hip{\mu - {\int \phi\ dP}}{f} + f_0 (1 - \pip 1 P)\\
= \pip lP -  \indc (\mu) + \hip{\mu}{f} - \pip {f}{ P} + f_0  - \pip {f_0}P\\
=\pip{l -f - f_0}{P} + \left(\hip{\mu}{f} -  \indc (\mu) \right) +f_0.
\end{multline}
The second equality is due to the reproducing property of RKHS.
The dual function is given by
$$
g(f,f_0) = \sup_{P, \mu} \pip{l -f - f_0}{P} + \left(\hip{\mu}{f} -  \indc (\mu) \right) +f_0.
$$
The first term  is bounded above by $0$ iff $l - f - f_0\in -K^*$.
By Lemma~\ref{thm:choquet}, this conic constraint
is equivalent to the constraint of \eqref{eq:dual_conjugate}, ${l(\xi)\leq  f_0 + f(\xi)   ,\ \forall \xi \in\mathcal X} $.

Finally, expressing the second term using convex conjugate 
$\suppc (f) = \sup_{\mu}\hip{\mu}{f} -  \indc (\mu)  $
concludes the derivation.
\end{proof}
Strong duality can potentially be adapted from the strong duality result of moment problem, e.g., \cite{shapiroDualityTheoryConic2001}.
However, we give a self-contained proof with only elementary mathematics that sheds light on the connection between the RKHS theory and 
distributionally robust optimization.
The proof is a generalization of the Euclidean space conic duality theorem (\cite{ben-talRobustOptimization2009a} Theorem~A.2.1) to infinite dimensions.
Figure~\ref{fig:separate} illustrates the idea of the proof.
\input{draft/strong_dual.Rmd}

Finally, we summarize the results above to prove the \kdro Generalized Duality Theorem~\ref{thm:dro_main_thm} 
\begin{proof}        
Theorem~\eqref{thm:dro_main_thm} is obtained by reformulating the inner moment problem in \eqref{eq:kdro} using the strong duality result in Proposition~\ref{thm:dual_conjugate}, i.e.,
\begin{multline}
        \begin{opt1}
                \min_{\var} \underset{P,\mu}{\sup}\bigg\{  \int l(\theta,\xi) \ dP(\xi) \colon\int{\phi }{\ dP} = \mu, P\in \mathcal P,\mu\in \mathcal C \bigg\}
        \end{opt1}
        \\
        =
        \min_{\var}\min_{f_0\in\mathbb R, f\in\rkhs}
        \bigg\{
        f_0 + \suppc (f): l(\var, \xi)\leq  f_0 + f(\xi)   ,\ \forall \xi \in\mathcal X
        \bigg\}
        ,
\end{multline}
which results in formulation~\eqref{eq:kdro_rc}.
\end{proof}

\subsection{Table~\ref{tbl:sets} deriving formulations for various choices of RKHS ambiguity set $\mathcal C$}
\input{draft/table_c.Rmd}
\input{draft/pf/prop2_shapes.md}


\subsection{Complementarity condition and proof}
\label{sec:complement}
\begin{lemma}[Complementarity condition]
        \label{thm:complement}
Let $P^*, f^*, f^*_0$ be a set of optimal primal-dual solutions of (P) and (D), then
\begin{equation}
        \pip{l-f^* - f_0^*}{P^*}=0 , \quad \suppc(f^*) =\pip{f^*}{ P^*   }
        \label{eq:complementarity_supp}
\end{equation}
\end{lemma}
If $\mathcal C=\{\mu\colon\|\mu - \mu_{\hat P}\|_\rkhs\leq \epsilon\}$, the second equality
implies
\begin{equation}                
        {\pip{\frac{f^*}{\hnorm{f^*}}}{(P^* - \hat{P})} = \mathrm{MMD}({P^* , \hat{P}})},
        \label{eq:cauchy_schwarz}
\end{equation}
which gives a second interpretation of the dual solution $f^*$ as a witness function.

It is well known that complementarity condition holds iff strong duality holds in the moment problem; cf. \cite{shapiroDualityTheoryConic2001}. The following is a straightforward proof. 
\begin{proof}
    Plug  $P^*, f^*, f^*_0$ into Lagrangian~\eqref{eq:lagrangian}, 
    $$
    \pip{l}{P^*} 
    \leq
    \pip{l -f - f_0}{P^*} + \suppc(f^*) +f_0^*
    \leq
    \suppc(f^*) +f_0^*.
    $$
    By strong duality, all inequalities above are equalities.
    Therefore, the first equality gives the condition $\suppc(f^*) =\pip{f^*}{ P^*   }$ while the second yields $\pip{l-f^* - f_0^*}{P^*}=0$.
\end{proof}

\subsection{Proof of Proposition~\ref{thm:rkhskiss} (Interpolation property)}
\begin{proof}
        Since $f_0^*, f^*$ is a solution to the inner moment problem of \kdro  \eqref{eq:dual_conjugate}, we have $l(\var, \xi)\leq  f^*_0 + f^*(\xi)   ,\ \forall \xi \in\mathcal X$ for any given \var. By the first equation in the complementarity condition~\eqref{eq:complementarity_supp}, we have $\pip{l-f^* - f_0^*}{P^*}=0 $. Hence the integrand must be zero ${P^*}$-a.e.
\end{proof}


\subsection{Corollary~\ref{thm:ipm_dual} IPM-DRO duality}
\label{sec:ipm}
We provide a derivation using a technique alternative to the proof of Proposition~\ref{thm:dual_conjugate}.
\begin{proof}
    We consider the Lagrangian
    \begin{multline}
        \mathcal L(P;\lambda) =
        \pip lP - \lambda (\ipm (P,\hat{P}) - \epsilon )
        \\
        =  \pip lP - \lambda \sup_{f\in \mathcal F}\int f d (P-\hat{P}) + \lambda\epsilon 
        \\
        =  \inf_{f\in \mathcal F} \pip {l - \lambda f }P + \lambda \int f d  \hat{P} + \lambda\epsilon \\
        \leq  \inf_{f\in \mathcal F} \sup_{\xi\in\domain} [{l(\xi) - \lambda f(\xi) }] +\frac{\lambda}N\sum_{i=1}^N f(\xi_i) + \lambda\epsilon
        .
    \end{multline}
    The second equality above is due to the dual representation of IPM. The last inequality is due to that the expectation is always dominated by the supremum.
    This results in the reformulation
    $$
    \min_{\var, \lambda\geq 0, f\in \mathcal F} \sup_{\xi\in\domain} [{l(\xi) - \lambda f(\xi) }] +\frac{\lambda}N\sum_{i=1}^N f(\xi_i) + \lambda\epsilon.
    $$
    By introducing the epigraphic variable $f_0$, we obtain the reformulation~\eqref{eq:ipm0}.
\end{proof}

\subsection{Corollary~\ref{thm:expect_cons} \kdro as stochastic optimization with expectation constraint}
Using the known relationship between semi-infinite constraint and expectation constraint (see, e.g., \cite[Theorem 1]{tadicRandomizedAlgorithmsSemiInfinite2006}), the SI constraint in \eqref{eq:kdro_rc} is equivalent to the expectation constraint in \eqref{eq:sa_exp_cons}.

%% file: draft/strong_dual.Rmd
\begin{proof}
We assume the dual optimal value of \eqref{eq:dual_conjugate} is finite $(D)<\infty$. 
Since the converse means that the dual problem is infeasible, which implies that the primal problem is unbounded. 
Due to the upper semicontinuity of $l$ in Assumption~\ref{asm:compact}, this can not happen on a compact $\domain$.

Let us consider the Hilbert space $\rset\times\rkhs\times\rset$ equipped with the inner product ${\langle, \rangle_\rset + \hip{}{} +\langle, \rangle_\rset}$.
We construct a cone in $\rset\times\rkhs\times\rset$
$$
A = \bigg\{ \left(\pip{1}{P}, \int\phi dP , \pip lP\right): P\in\mathrm{co}(\mathcal P) \bigg\},
$$
where $\mathrm{co}$ again denotes conic hull.

Let $t=(P)$ denote the optimal primal value. $\forall \epsilon>0$, we construct the set
$$
B_\epsilon =\bigg\{ \left(1, \mu , t+\epsilon \right): \mu\in
 \mathcal C
\bigg\},
$$
which is a closed convex set with non-empty relative interior by Assumption~\ref{asm:compact} (i.e., Slater condition is satisfied).

It is straightforward to verify that those two sets do not intersect.
Suppose $x=(x_1, x_2, x_3) \in A\cap B_\epsilon$, this means $\exists \mu',P'$ such that 
$x_1 = 1 = \pip 1{P'}, x_2 =\mu' = \int\phi d{P'}$,
 i.e., $ \mu',P'$ is a primal feasible solution. Then the third coordinate of $x$ satisfies $x_3 = \pip l{P'} \leq (P) < t+\epsilon = x_3$, which is impossible. 
Hence, $A\cap B_\epsilon = \emptyset$.

In the rest of the proof, we will show that, $\forall \epsilon>0$, the dual optimal value $(D)$ satisfies
$$
(D) \leq (P) + \epsilon.
$$
Combining this with weak duality $(D) \geq (P)$ will result in strong duality. We now justify this inequality.

By the separation theorem, 
(see, e.g., \cite{barvinokCourseConvexity2002} Theorem~III.3.2, 3.4), there exists a closed hyperplane that strictly separates $A$ and $B_\epsilon$. 
The separation is strict because $t+\epsilon > t = \pip lP$.
By the Riesz representation theorem, $\exists (f_0, f, \tau)\in\rset\times\rkhs\times\rset, s\in\rset$, such that
\begin{align*}
f_0 + \hip{f}{\mu} +\tau( t+\epsilon ) < s,\\
f_0 \pip{1}{P} + \hip{f}{\int\phi dP} +\tau\pip lP > s.
\end{align*}
Plugging in $P=0$, we obtain $s < 0$. Since $P$ lives in a cone, the left-hand side of the second inequality must be non-negative. Otherwise, we can scale $P$ so that the separation will fail. 
In summary, we have
\begin{equation}
    \begin{aligned}
f_0 + \hip{f}{\mu} +\tau( t+\epsilon ) & <  0,\forall \mu\in\mathcal C ,\\
f_0 \pip{1}{P} + \hip{f}{\int\phi dP} +\tau\pip lP & \geq  0, \forall P\in\mathrm{co}(\mathcal P) .
\end{aligned}
\label{eq:strong_pf1}
\end{equation}

By Assumption~\ref{asm:compact} (Slater condition),
the primal problem has a non-empty solution set. 
Because $l$ is proper and upper semi-continuous 
and the feasible solution set for the optimization problem is compact (see Section~\ref{sec:lemmas})
,
the primal optimum is attained by the extreme value theorem.
Suppose $P^*$ is a primal optimal solution, from the second inequality of \eqref{eq:strong_pf1},
\begin{align*}
f_0 \pip{1}{P^*} + \hip{f}{\int\phi dP^*} +\tau\pip l{P^*} 
=f_0  + \hip{f}{\mu_ {P^*}} +\tau t \geq 0.
\end{align*}
Using this and the first inequality of \eqref{eq:strong_pf1}, we obtain $\tau < 0$. 
Without loss of generality, we let $\tau=-1$.

From the second inequality of \eqref{eq:strong_pf1}, we have
$$
f_0 \pip{1}{P} + \hip{f}{\int\phi dP}  - \pip lP =  \pip{f_0 + f -l}P
\geq 0, \forall P\in\mathrm{co}(\mathcal P).\\
$$
This tells us that $f_0, f$ is a feasible dual solution because it satisfies the semi-infinite constraint in \eqref{eq:dual_conjugate}.

By the first inequality of \eqref{eq:strong_pf1},
$$
f_0 + \sup_{\mu\in\mathcal C}\hip{f}{\mu} \leq   t+\epsilon  , \forall\epsilon>0,
$$
where the left-hand side is precisely the dual objective in \eqref{eq:dual_conjugate}.
This implies $(D) \leq (P) + \epsilon.$ By weak duality, $(D) \geq (P)$. Therefore, strong duality holds.
\end{proof}
This proof gives us the third interpretation of the dual variables $f_0, f$ --- they define a separating hyperplane of $A$ and $B_\epsilon$.
\begin{remark}
    From the proof, we see that the \emph{Slater condition} in Assumption~\ref{asm:compact} is stronger than needed be.
    If $\mathcal C$ is singleton, we can still find a convex neighborhood $W_\epsilon$ of the singleton $B_\epsilon$ since $\rset\times\rkhs\times\rset$ is locally convex. Then $W_\epsilon$ and $A$ can still be strictly separated using the same technique in the proof. Hence strong duality still holds when $\mathcal C$ is a singleton.
\end{remark}

%% file: draft/table_c.Rmd
\begin{table*}
        \caption{Robust counterpart formulations of \kdro. }
        \label{tbl:sets}
        \centering
        \begin{tabular}{ll}
          \toprule
        RKHS ambiguity set $\mathcal C$ & Robust counterpart formulation\\
        \midrule
          norm-ball $\mathcal C=\{\mu\colon\|\mu - \mu_{\hat P}\|_\rkhs\leq \epsilon\}$ &
          $
          {f_0+\frac1N\sum_{i=1}^N f(\xi_i)+\epsilon\hnorm{f} }
          $
          \\

          \midrule
          convex hull $\mathcal C=\mathrm{conv}\{\mathcal C_1,\dots,\mathcal C_N\}$& 
          ${f_0+\max_i \delta^*_{\mathcal C_i} (f)}$
          \\
          (same under closure $\mathrm{cl}\mathrm{conv}\{\mathcal C \}$
          )& 

          \\
          \addlinespace[.5em]
          example: polytope
          & 
          ${f_0+\max_i f(\xi_i)}$
          \\
          $\mathcal C=\mathrm{conv}\{\phi(\xi_1),\dots,\phi(\xi_N)\}$
          &
          (equivalent to SVMs/scenario opt. \cite{calafioreScenarioApproachRobust2006})
          \\
          \midrule
          Minkowski sum           $\sum_{i=1}^N C_i$ & $f_0+ \sum_{i=1}^N\delta^*_{\mathcal C_i} (f)$\\
          \addlinespace[.5em] 
          example: $\mathcal C = \mathcal C_1+\mathcal C_2$ 
          &
          ${ f_0+\max_i f(\xi_i)+\epsilon\hnorm{f}  }$
          \\
          $\mathcal C_1=\{\mu\colon\|\mu \|_\rkhs\leq \epsilon\}$\\$\mathcal C_2 = \mathrm{conv}\{\phi(\xi_1),\dots,\phi(\xi_N)\}$
          \\
          \midrule
          affine combination 
          & 
          ${f_0+\sum_{i=1}^N\alpha_i \delta^*_{\mathcal C_i} (f)}$
          \\
          $\mathcal C=\sum_{i=1}^N\alpha_i \mathcal C_i, \sum_{i=1}^N\alpha_i=1$
          \\
          example: data contamination
          &
          ${f_0+\frac\alpha N\sum_{i=1}^N f(\xi_i)} + (1-\alpha) \delta^*_{\mathcal C_Q}(f)$
          \\          
          $\mathcal C=\{\alpha \mu_{\hat{P}} + (1-\alpha) \mu_{Q}: \mu_{Q}\in \mathcal C_Q\}$
          \\
          \midrule
          Intersection    $\mathcal C =\cap_{i=1}^N\mathcal C_i$ 
          & $f_0+ \sum_{i=1}^N\delta^*_{\mathcal C_i} (f_i),\ \sum_{i=1}^Nf_i = f$
          \\
          \midrule
          multiple kernels $\mathcal C_i\subseteq\rkhs_i$  
          &
          $
          {f_0+ \sum_{i=1}^N\delta^*_{\mathcal C_i} (f_i)}$ where $f_i\in\rkhs$
          \\
          \addlinespace[.5em] 
          example: 
          $
          \mathcal C_i=\{\mu\colon\|\mu - \mu_{\hat P}\|_\rkhs\leq \epsilon_i\}
          $
          &
          ${f_0+\frac1N\sum_{i=1}^N\sum_{i=j}^N f_i(\xi_j)+\epsilon \sum_{i=1}^N\|{f_i}\|_{\rkhs_i}}$
          \\
          \midrule
          singleton 
          $\mathcal C = \left\{\kme{\frac{1}{N}}{\xi_i}\right\}$
          &
           $f_0+\frac1N\sum_{i=1}^N f(\xi_i)$ (equivalent to ERM/SAA)
           \\
        \midrule
        entire RKHS $\mathcal C=\rkhs$ &
        $f_0+$ $\delta_0(f)$\\
        &
        (equivalent to worst-case RO \cite{ben-talRobustOptimization2009a})
        \\
          \bottomrule
        \end{tabular}
\end{table*}

%% file: draft/pf/prop2_shapes.md
We now derive the formulations of support functions for various RKHS ambiguity sets in Table~\ref{tbl:sets}.
\paragraph{(RKHS norm-ball) }
Let us consider the ambiguity set of
$\mathcal C=\{\mu\colon\|\mu - \hat\mu\|_\rkhs\leq \epsilon\}$, where ${\hat P}=\sum_{i=1}^N{\frac1N}\delta_{\xi_i}$.
The support function is given by
$$
\suppc (f) 
= \sup_{\mu\in\mathcal C} \hip{f}{\mu} =  \hip{f}{\hat\mu} + \sup_{\|\mu - \hat\mu\|_\rkhs\leq \epsilon} \hip{f}{\mu-\hat\mu} = \hip{f}{\hat\mu} + \epsilon\hnorm{f}
$$
where the last equality is by the Cauchy-Schwarz inequality, or alternatively by the self-duality of Hilbert norms. 
(Note we assume 
there exists some $\mu\in\rkhs$ such that $\|\mu - \hat\mu\|_\rkhs=\epsilon$.)

\paragraph{(Polytope, convex hull of ambiguity set)}
The result for convex hull follows from standard support function calculus.
If the ambiguity set $\mathcal C$ is described by the polytope $\mathrm{conv}\{ \phi(\xi_1),\dots,\phi(\xi_N)\}$,
then $\suppc (f) = \max_{1\leq i\leq N} f(\xi_i)$.
Furthermore, the support function value remains the same under closure operation
\footnote{The convex hull can be replaced with its closure $\mathrm{cl conv}(\cdot)$. Note convex hulls in infinite-dimensional spaces are not automatically closed; cf. Krein-Milman theorem.}
.
The equivalence to the scenario approach in \cite{calafioreScenarioApproachRobust2006} can be seen by noticing that
$
\max_{1\leq i\leq N} l(\xi_i) \leq f_0+  \max_{1\leq i\leq N} f(\xi_i)  .
$
If $\rkhs$ is universal, then there exists $f_0, f$ such that the equality is attained.
\paragraph{(Minkowski sum, affine combination, intersection)} Those cases follow directly from the support function calculus; cf. \cite{ben-talDerivingRobustCounterparts2015}.

\paragraph{(\kdro with multiple kernels)} Let us consider multiple ambiguity sets from different RKHSs.
Suppose $\rkhs_1,\dots, \rkhs_{N_h}$ are RKHSs associated with feature maps $\phi_1,\dots, \phi_{N_h}$.
Let $\mathcal{C}_1,\dots, \mathcal{C}_{N_h}$ be the ambiguity sets in the respective RKHSs.
\kdro formulation with multiple kernels is given By
\begin{equation}
                \min_{\var} \underset{P,\mu}{\sup}\bigg\{{  \int l(\theta,\xi) \ dP(\xi) }\colon
                {\int{\phi_i }{\ dP} = \mu_i, P\in \mathcal P,\mu_i\in \mathcal C_i}, i=1\dots N_h \bigg\},
\end{equation}
Using the same proof as Proposition~\ref{thm:dual_conjugate}, we have the \kdro reformulation
\begin{equation}
    \begin{opt1}
            \MINST{ \var,f_0\in\mathbb R, f_i\in\rkhs_i}{\  f_0 + \sum_{i=1}^N{\delta^*_{\mathcal{C}_i}} (f_i)  }
            {
                l(\var, \xi)\leq  f_0 + \sum_{i=1}^Nf_i(\xi)   ,\ \forall \xi \in\mathcal X    
            } .
    \end{opt1}
    \end{equation}
Hence we obtain the formulation in Table~\ref{tbl:sets}.
\paragraph{(Singleton ambiguity set $\mathcal C = \left\{\kme{\frac{1}{N}}{\xi_i}\right\}$)}
By the reproducing property, the support function of the singleton ambiguity set is given by $\suppc (f) = \frac1N\sum_{i=1}^N f(\xi_i)$.
\paragraph{(If $\mathcal C=\rkhs$, reduction to classical RO)}
$\delta^*_\rkhs (f) \neq \infty$ iff $f=0$. Then \eqref{eq:kdro_rc} is reduced to
\begin{equation}
    \begin{opt1}
            \MINST{ \var, f_0\in\rset}{\  f_0 }{l(\var, \xi)\leq f_0 ,\ \forall \xi \in\mathcal X} 
    \end{opt1}
\end{equation}
which is the epigraphic form of the worst-case RO.
\footnote[2]{Note that $\mathcal C=\rkhs$ is no longer closed. However, the resulting ambiguity set becomes $\mathcal P$, which is still compact if \domain is compact.}

%% file: draft/randfeat.md
Common ways to parametrize an RKHS function include the representer theorem as well as approximations such as the random Fourier features~\cite{rahimiRandomFeaturesLargeScale2008}.
Recall that an RKHS function can be approximated by the finite feature expansion
$$
f(\xi) \approx \hat{f}(\xi)= w^\top \hat\phi(\xi), \ \ k(x, x')\approx
\sum_{i=1}^N \hat\phi_i(x)\hat\phi_i(x')
$$

where $\{\hat\phi_i(x)\}_{i=1}^N$ are the random features, e.g., random Fourier features
$
\hat\phi_i(x) = \cos(w_i x+b_i), w_i\sim \mathrm{N}(0, \sigma^2), b_i\sim \mathrm{Uniform}[0,2\pi].
$
If $x$ is a vector, then $w_i \sim \mathcal{N}(0, I \sigma^2)$, and $w_i  x$ is the dot product.
See, e.g., \cite{rahimiRandomFeaturesLargeScale2008}, for more properties.

One strength of the Generalized Duality Theorem~\ref{thm:dro_main_thm} is that it does not require the knowledge of the RKHS that the loss $l$ lives in, which is typically not available in non-kernelized models.
This enables us to use approximate features for commonly used RKHSs, e.g., random Fourier feature.
This is a strength of our \kdro theory.

Note program~\eqref{eq:scenario} is a convex optimization problem with the random feature parametrization.

\subsection{Distributionally robust version of representer theorem}
In program~\eqref{eq:scenario},
we may parametrize the RKHS function by 
$f(\cdot) = \sum_{i=1}^{N}{\beta_i}{k(\xi_i,\cdot)} + \sum_{j=1}^M{\gamma_j}{k(\zeta_j,\cdot)}$, 
$\hnorm{f} = \sqrt{\alpha^\top K \alpha}$, where $\alpha = (\beta_1, \ldots, \beta_N, \gamma_1, \ldots, \gamma_M)^\top, K = [k(\eta_i,\eta_j)], \eta = (\xi_1, \ldots, \xi_N, \zeta_1, \ldots, \zeta_M)^\top$.
We justify this parametrization by
the following DRO version of the RKHS representer theorem~\cite{SchHerSmo01}.

The intuition of the following result is to restrict \kdro to a smaller ambiguity set of distributions supported on $\{\zeta_i\}_{i=1}^M$ (i.e., replace $P\in\mathcal P$ by $P\in\mathcal P_M$, an inner approximation depending on $M$).
In this setting,
the ambiguity set only contains only distributions supported on (a {subset} of) $\zeta_i$.
Then it suffices to parametrize $f$ in \eqref{eq:scenario} by $f(\cdot) = \sum_{i=1}^{N}{\beta_i}{k(\xi_i,\cdot)} + \sum_{j=1}^M{\gamma_j}{k(\zeta_j,\cdot)}$.

\begin{lemma}[Robust representer]
\label{thm:representer_weak}
Given data $\{\xi_i\}_{i=1}^N$ and the ambiguity set chosen to be a set of embeddings with the form $\sum_{j=1}^M\alpha_j\phi(\zeta_j)$, for some $0\leq \alpha_j\leq 1, \sum_{j=1}^M\alpha_j=1$, and within the
RKHS norm-ball
$\mathcal C=\{\mu\colon\|\mu - \mu_{\hat{P}}\|_\rkhs\leq \epsilon\}$.
Then,
it suffices to consider the RKHS function of the form
$
  f(\cdot) = \sum_{i=1}^{N}{\beta_i}{k(\xi_i,\cdot)} + \sum_{j=1}^M{\gamma_j}{k(\zeta_j,\cdot)}
$
  for some $\beta_i,\gamma_j \in\mathbb R, i=1...N, j=1...M$.
\end{lemma}
Lemma~\ref{thm:representer_weak}
states that the expansion points of the RKHS representer in \eqref{eq:scenario} are exactly the support of the empirical and worst-case distributions.
It extends the classical RKHS representer theorem~\cite{SchHerSmo01}, which uses only the empirical samples as expansion points. The implication is that, to be distributionally robust, we should choose the representers as in Lemma~\ref{thm:representer_weak} instead of only using empirical samples.
Below is a proof that is similar to the original representer theorem.
\begin{proof}
    In \eqref{eq:scenario} we consider $f = f_s + f_\perp$, where $ f_s= \sum_{i=1}^{N}{\beta_i}{k(\xi_i,\cdot)} + \sum_{j=1}^M{\gamma_j}{k(\zeta_j,\cdot)}$ belongs to a subspace of the \rkhs and $f_\perp$ its complement.
    Plug in $f = f_s + f_\perp$ to \eqref{eq:scenario} and note the orthogonality, we obtained
    \begin{equation}
                \begin{opt1}
                \MINST{\var, f_s,f_\perp, f_0}
                {f_0+\frac1N\sum_{i=1}^N f_s(\xi_i) + \epsilon(\hnorm{f_s}+\hnorm{f_\perp} )}
                {l(\var,\zeta_i) \leq f_s(\zeta_j)  + f_0 ,\  j=1\dots M.}
                \end{opt1}
    \end{equation}
    It suffices to choose $f_\perp=0$ in this optimization problem. Hence the conclusion follows.
\end{proof}

\begin{remark}
        Note the existence of a worst case distribution in more general settings is not yet proven. The discussion here is restricted to the setting of \eqref{eq:scenario}.
\end{remark}

%% file: draft/certify.Rmd
We carry out additional numerical experiments to study \kdro.

\subsection{Testing other variants of \kdro}
We empirically test the following proposed variants of \kdro.
\begin{itemize}
	\item Relaxed \kdro formulation (\kdro-relaxed) with constraint hold for only the empirical samples, i.e., 
	${l(\var,\xi_i) \leq f_0 + f(\xi_i),
	\ i=1\dots N.}$
	\item Unconstrained \kdro using Kernel CVaR in Example~\ref{ex:kcvar}
\end{itemize}
\begin{figure}[h!]
	\centering
	\includegraphics[width=.55\columnwidth]{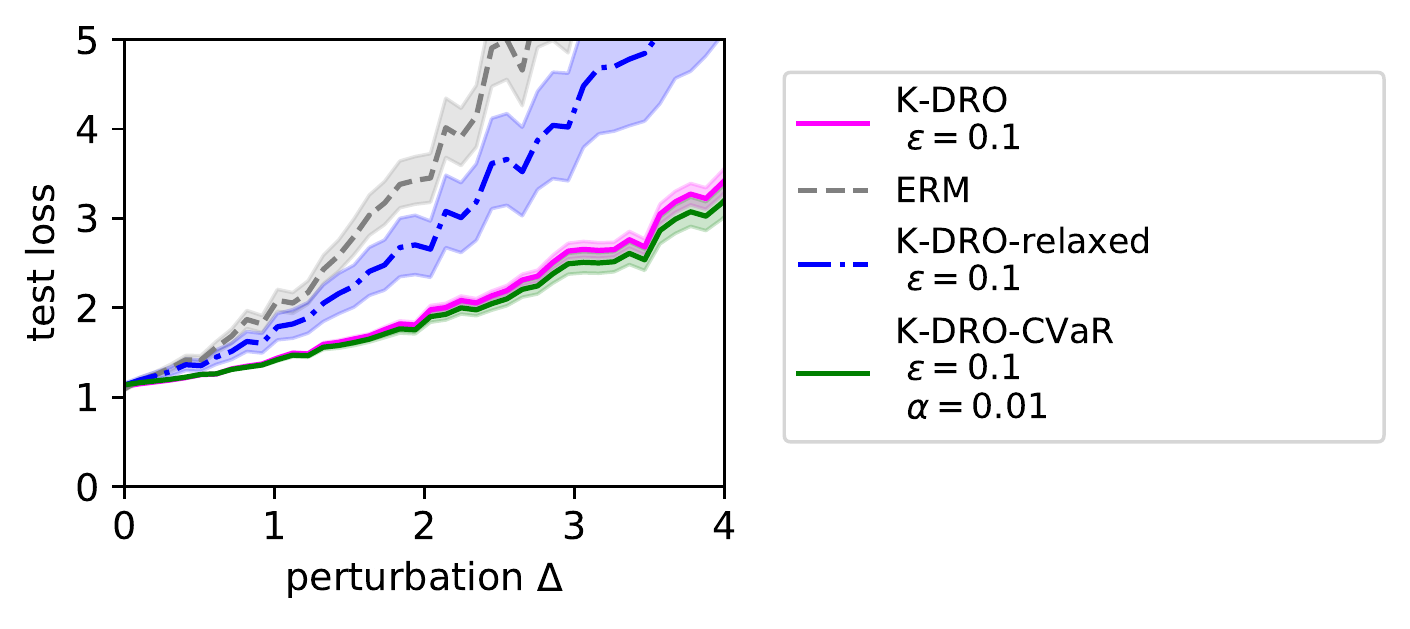}
	\caption{Comparing \kdro-relaxed, \kdro-KCVaR, ERM, and regular \kdro. y-axis limit is adjusted to show the plot. All error bars are in standard error.
			}
	\label{fig:no_certify}
\end{figure}

We compare \kdro-relaxed with the ERM as well as the regular \kdro. 
Compared with ERM, \kdro-relaxed still possesses moderate robustness. In this case, we effectively proposed a way to apply RKHS regularization to general optimization problems, not limited to kernelized models.
Hence, it may be used in practice as a finite-sample approximation to \kdro.

We then test the \kdro using the unconstrained objective given by Kernel CVaR. We observe no significant difference in performance between \kdro-KCVaR (with small chance constraint level $\alpha$.) and regular \kdro~\eqref{eq:scenario}.

\subsection{Analyzing the generalization behavior}
\begin{figure}[htb!]
	\centering
	\includegraphics[height=4cm]{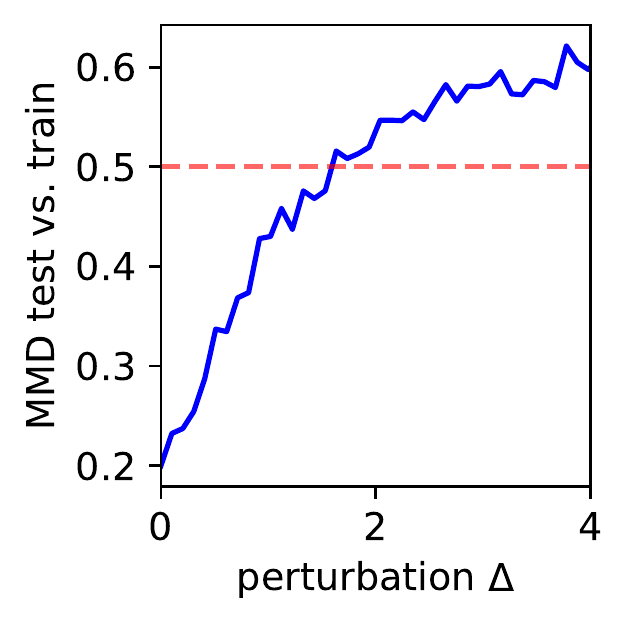}
  \hspace{1cm}
	\includegraphics[height=4cm]{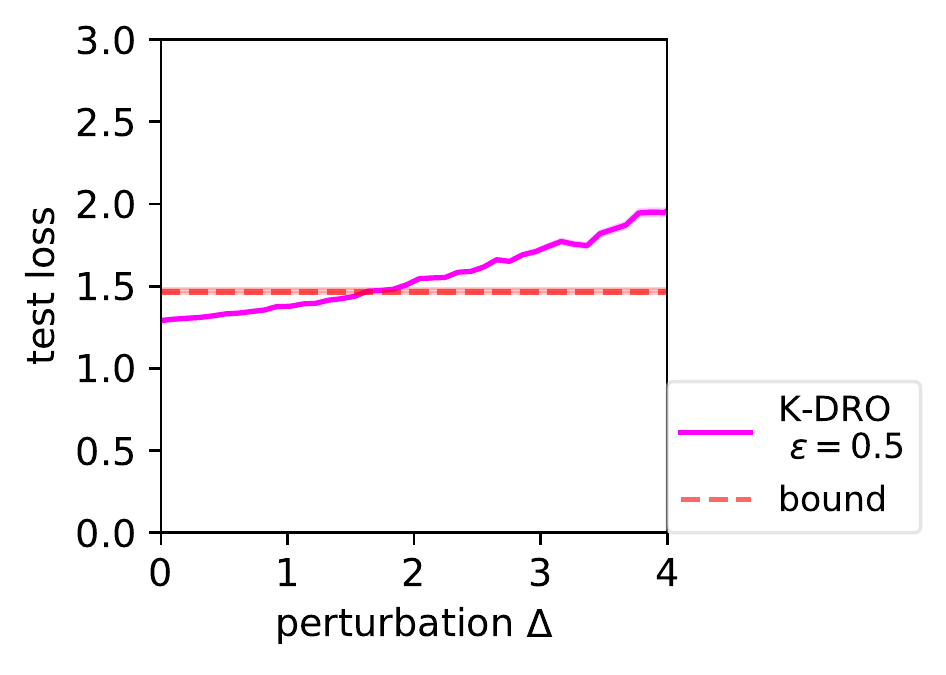}
	\caption{
		(Left) MMD estimator between the empirical samples and test samples. The level $\mathrm{MMD}=0.5$ is marked in red.
		(Right) Loss compared to the generalization bound. As the test data falls outside the robustification level $\epsilon$, the loss starts to exceed the generalization bound (red) $
		{ f_0+\frac1N\sum_{i=1}^N f(\xi_i)+\epsilon\hnorm{f}}$.
		}
		\label{fig:exp_mmd}
\end{figure}

An insight can be obtained by observing the plot of the MMD estimator between the training and test data in Figure~\ref{fig:exp_mmd} (left). As \kdro with $\epsilon=0.5$ robustified against perturbation less than the level $\mathrm{MMD}=0.5$, we see this threshold was exceeded as we increase the perturbation in test data. Meanwhile, this is the same time ($\Delta\approx 1.5$) where \kdro solutions start to exceed the generalization bound~ 
$
{\int l \ dP \leq f_0+\frac1N\sum_{i=1}^N f(\xi_i)+\epsilon\hnorm{f}}
$, see Figure~\ref{fig:exp_mmd} (right). 
This empirically validates our theoretical results for robustification.


\subsection{Miscellaneous details for experimental set-up}
\label{sec:misc_exp}
\paragraph*{Robust least squares example.}
Our experiments are implemented in Python.
The convex optimization problems are solved using ECOS or MOSEK interfaced with CVXPY. 
In the experiments, we chose the bandwidth for the Gaussian kernel using the medium heuristic~\cite{grettonKernelTwosampleTest2012}. $\epsilon$ in this paper are fixed to constants below $2$ for Gaussian kernels. Choosing $\epsilon$ can be further motivated by kernel statistical tests~\cite{grettonKernelTwosampleTest2012} and is left for future work.

\paragraph{Sampled $\zeta_j$}
In applying \kdro using ~\eqref{eq:scenario},
we may obtain $\zeta_j$ by simply sampling in \domain.
$\{ \zeta_j \}_i$ need not be real data, e.g.,
in stochastic control, they can be a grid of system states;
in learning, they can be
synthetic samples such as convex combinations of data $\zeta_j = \sum_{i=1}^N a_{ij} \xi_i, a_{\cdot j}\in S_N$ (simplex), or
perturbations $\zeta_j = \hat{\xi}_i+\Delta_i$ where $\Delta_i$ can be a small perturbation,
or they can be obtained by domain knowledge of the specific application.
In the setting of supervised machine learning, there is a difference between this paper's approach of sampling $\zeta_j$ and commonly used data-augmentation techniques: $\zeta_j$ need not have the correct labels or targets.
Directly training on them may have unforeseen consequences.
For example, in the robust least squares experiment, we sampled the support $\zeta_i$ uniformly random from $[-1, 1]$.

\paragraph*{Robust learning under adversarial perturbation example.}
For the MNIST robust classification example, we used a neural network with two hidden layers with $64$ units each.
For the training of ERM and PGD, we used the ADAM optimization routine implemented in the PyTorch library.
In Step~\ref{step:feature} of \sfg, we used random Fourier features~\cite{rahimiRandomFeaturesLargeScale2008} with $500$ features.
In Step~\ref{step:sa} of \sfg, we used the SA routine from CSA algorithm~\cite{lanAlgorithmsStochasticOptimization2020a}.
While other SA routines can be used, we prefer the simplicity of CSA in that it does not use a dual variable.
We set the threshold and step-size of the CSA algorithm~\cite{lanAlgorithmsStochasticOptimization2020a} to decay at the rate of $\frac{1}{\sqrt{k}}$ as suggested in that paper.
We did not attempt further adaptive tuning of the step-sizes or the proposing distribution for $\zeta$ (we generate $3000$ samples uniformly in Step~\ref{step:sample} of \sfg), which may further improve the performance. Parameter (weights of the neural nets) averaging is used for training all models.
In the visualization of the predictions in Figure~\ref{fig:mnist_viz}, we perturbed the images by the PGD method~\cite{madryDeepLearningModels2019,madryAdversarialRobustnessTheory} based on the ERM loss and linear model. SFG-DRO does not have the knowledge of the perturbation method.

%% file: draft/pf/gap.Rmd
\subsection{Computing worst-case distributions}
\label{sec:recover}
\begin{figure}[t]
                \centering
                \subfloat[Support sampling\label{fig:perturb_0_sample}]{
                \includegraphics[width=.3\columnwidth,valign=c]{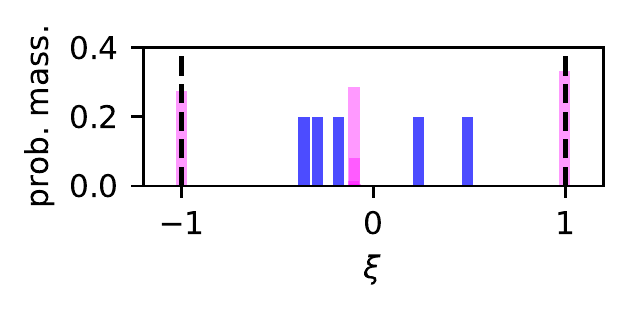}
                \includegraphics[width=.05\columnwidth, valign=c]{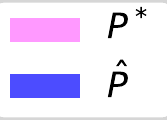}
                }
                \subfloat[Perturbation\label{fig:perturb_0_opt}]{
                \includegraphics[width=.3\columnwidth,valign=c]{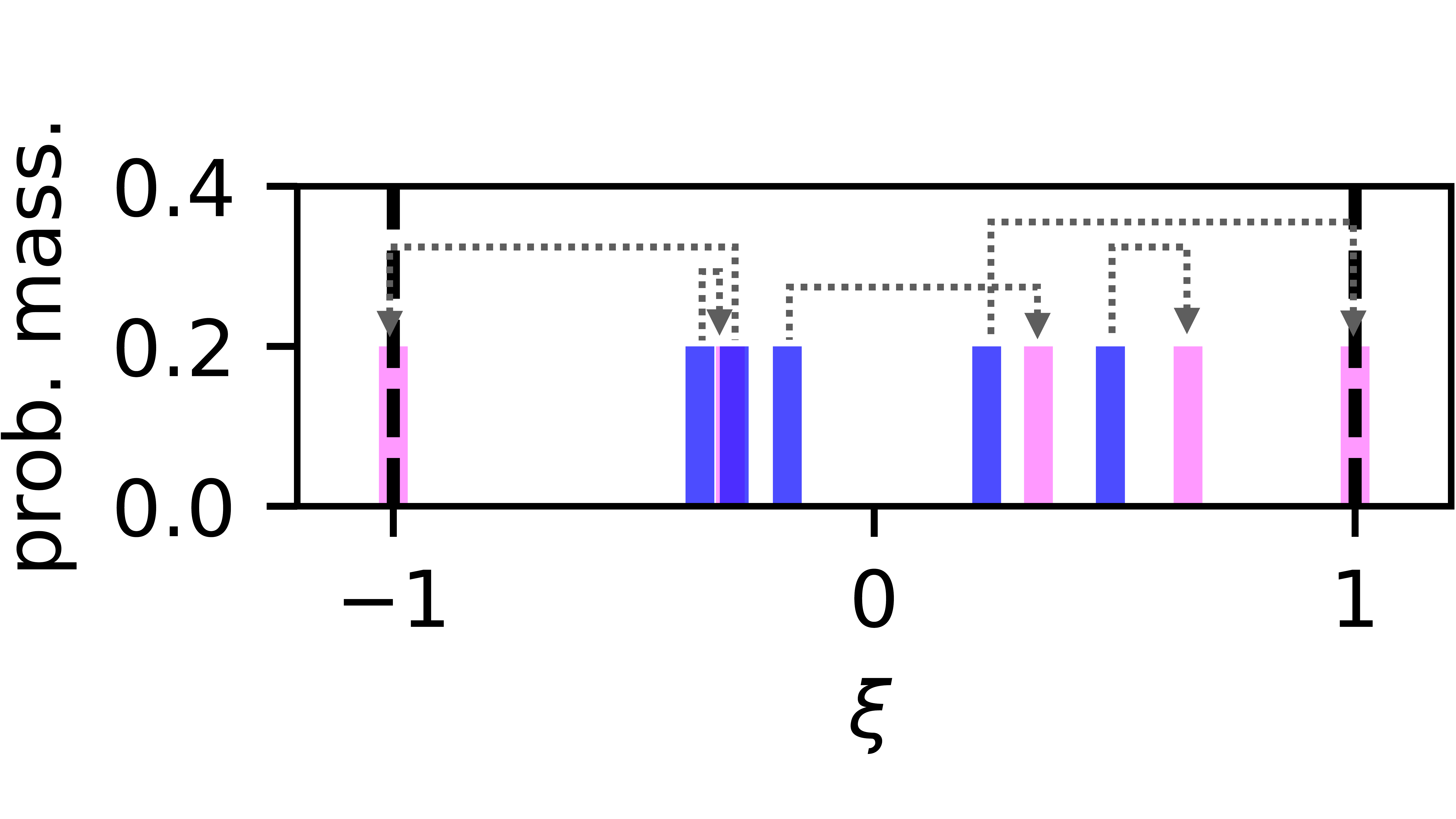}
                }
                \caption{Computing the worst-case distribution $P^*$ using
                        \textbf{(a)}: \eqref{eq:cdc} samples possible support $\zeta_j$ then optimizes w.r.t. weights $\alpha$.
                \textbf{(b)}: \eqref{eq:perturb} moves the empirical samples directly.
                }
                \label{fig:perturb_0}
                \vspace{-5mm}
        \end{figure}

We have proposed \kdro for making the decision \var via reformulation~\eqref{eq:kdro_rc}. In practice, it is often useful to find the worst-case distribution $P^*$ (e.g., to study adversarial examples).
We now propose two practical methods to compute $P^*$
for a given \var,
based on \emph{support sampling} and \emph{perturbation}, respectively.
We illustrate the ideas in Figure~\ref{fig:perturb_0}. 
\vspace{-0.3cm}
\paragraph{Support sampling.}
We consider the moment problem \eqref{eq:moment_primal} where the distribution is restricted to discrete distributions supported on some sampled support $\{\zeta_j\}_{j=1}^M\subseteq\domain$.
\footnote{
\label{foot:support}
Note the sampled support $\{\zeta_j\}_{j=1}^M $ need not be real data; they are only the candidates for the worst-case support.
The purpose is to make the the semi-infinite constraint approximately satisfied.
See the appendix for more details.
}
For any given \var,
\begin{equation}        
    \mxz_{\alpha\in S_{M}} \sum_{j=1}^M\alpha_i l(\var, \zeta_j)\quad
    \sjt \ 
    \bigg\Vert{\sum_{j=1}^M\alpha_i \phi(\zeta_j) - \frac1N\sum_{i=1}^N \phi(\xi_i)}\bigg\Vert_\rkhs
    \leq\epsilon. \\
    \label{eq:cdc}
\end{equation}
\eqref{eq:cdc} can be written as a quadratically constrained program with linear objective,
which admits a (strong) semidefinite program dual
via what is historically known as the S-lemma \cite{polikSurveySLemma2007} (cf. appendix).
Alternatively, \eqref{eq:cdc} can be directly handled by convex solvers for a given \var.
Note this approach was previously used in solving the problem of moments in~\cite{zhuWorstCaseRiskQuantification2020}.
\vspace{-0.3cm}
\paragraph{Perturbation.}
Alternatively, we search for worst-case distributions that are \emph{perturbations} of the empirical distribution. 
Let $d_i\in\domain$ be some perturbation vector, given \var,
\begin{equation}        
        \begin{aligned}
        &\mxz_{\substack{d_i, i=1\dots N,\\ {\xi}_i + d_i\in\domain}} 
        & &
        \frac1N\sum_{i=1}^N l(\var, {\xi}_i + d_i)
        & \sjt
         & &
         \bigg\Vert{\frac1N\sum_{i=1}^N 
         \left(
         \phi({\xi}_i + d_i) 
         - 
         \phi({\xi}_i)
         \right)
         }\bigg\Vert_\rkhs
         \leq\epsilon.
        \end{aligned}
        \label{eq:perturb}
\end{equation}
Compared with \eqref{eq:cdc}, \eqref{eq:perturb} directly searches for the support of the worst-case distribution. It can be interpreted as transporting the probability mass from empirical samples $\xi_i$ to form the worst-case distribution.
Depending on the kernel used, \eqref{eq:perturb} may become a nonlinear program. However, its feasibility is guaranteed since it can always be initialized with a feasible solution $d_i=0$.

We now empirically examine the \emph{support sampling} method~\eqref{eq:cdc} and \emph{perturbation} method~\eqref{eq:perturb} to recover the worst-case distribution.
Since both programs \eqref{eq:cdc} and \eqref{eq:perturb} 
search for the worst-case distribution within a subset of all distributions,
their optimal values lower-bound the true worst-case risk (P) in \eqref{eq:moment_primal}, i.e., with finite samples, they are optimistic bound.

Under the experimental setting as in Figure~\ref{fig:lsq_interpret}, we ran \kdro with fewer empirical samples ($N=5$). After we obtain the \kdro solution $\var^*$, we plug it into \eqref{eq:cdc} and \eqref{eq:perturb}, respectively, to compute the worst-case distribution $P^*$. Figure~\ref{fig:perturb_0} plots the results. Note \eqref{eq:cdc} is a convex optimization problem, while \eqref{eq:perturb} results in a nonlinear program (with Gaussian kernel). Nonetheless, we solve it with an always-feasible initialization $d_i=0$.


\subsection{SDP dual via S-lemma}
\label{sec:sdp}
We consider a discretized version of the primal moment problem in \eqref{eq:cdc} where the distribution is constrained to be a discrete distribution. We rewrite \eqref{eq:cdc} as a quadratically constrained program using the plug-in estimator of MMD,
\begin{equation*}        
    \begin{aligned}
    &\mxz_{\alpha} & &  \sum_{i=1}^M\alpha_i l(\zeta_i) \\
    &\sjt& &\alpha^\top K_z\alpha - 2\frac1N \alpha^\top K_{zx} \mathbf{1} + \frac{1}{N^2}\mathbf{1}^\top K_x\mathbf{1} \leq\epsilon^2  \\
    &&& \sum_{i=1}^M\alpha_i =1, \alpha_i\geq0, i=1\dots M.
    \end{aligned}
\end{equation*}
This is a quadratically constrained linear objective convex optimization problem, 
where the Gram matrix $K_z$ almost always has exponentially decaying eigenvalues.
By applying S-lemma~\cite{polikSurveySLemma2007}, this program can be reformulated as the following SDP,
\begin{equation}
    \begin{opt}
            \MIN{\lambda\geq 0,x,y\geq0, t}{t}
            \SJT{
                    \begin{bmatrix}
                    \lambda P & -\lambda q - \frac12 (l + x\cdot \mathbbm 1 + y)\\
                    (-\lambda q - \frac12 (l + x\cdot \mathbbm 1 + y))^\top & t - \lambda\epsilon^2+x+\lambda r
                    \end{bmatrix}
                    \geq
                    0,
                    }
    \end{opt}
    \label{eq:sdp}
\end{equation}
where 
$
P:=K_z,\ q:= \frac1N  \mathbf{1}^\top K_{zx},\ r:=\frac{1}{N^2}\mathbf{1}^\top K_x\mathbf{1},
$
and $K_z= [k(\zeta_i, \zeta_j)]_{ij}, K_{zx}= [k(\zeta_i, \hat\xi_j)]_{ij}, K_x= [k(\hat\xi_i, \hat\xi_j)]_{ij}, l = [l(\zeta_1),\dots,l(\zeta_M)]^\top$.


%% file: draft/compact.Rmd
We establish a few technical results that are used in the proofs.
\subsection{Reducing conic constraint to infinite constraint}
To derive the semi-infinite constraint in \eqref{eq:dual_conjugate}, we need a standard result from the literature of the moment problem.
We give a self-contained proof below.
\begin{lemma}
        \label{thm:choquet}
        Let $K^*$ be the dual cone to the probability simplex $\mathcal P$.
        The conic constraint ${l - f - f_0\in -K^*}$ is equivalent to
        \begin{equation}
                \semiinfcons .
        \label{eq:conic_dual_choquet}
        \end{equation}
\end{lemma}
\label{sec:choquet_reduction}
\begin{proof}
``$\implies$'':  
Let us consider the set of all Dirac measures on \domain, $\mathcal D:=\{ \delta_{\xi}\colon \xi \in \domain\}$. For any $\xi\in\domain$, we have
$$
l(\xi) - f_0 - f(\xi) = \int l - f_0 - f d \delta_{\xi} \leq 0.
$$
Hence sufficiency.

``$\impliedby$'': Suppose there exists $P' \in\mathrm{co}(\mathcal P)$ such that $\int l - f_0 - f d P' > 0$. Without loss of generality, we assume $P'\in\mathcal P$, or we can normalize it to be a probability measure. Then,
$$
0 < \int l - f_0 - f d P' \leq \sup_{\xi\in\domain} l(\xi) - f_0 - f(\xi) \leq 0.
$$
The second inequality is due to that expectation is always less than or equal to the supremum. The last inequality holds because $l$ is u.s.c. This double inequality is impossible, hence $ l - f_0 - f \in -K^*$.
\end{proof}
Note an extension of this result to generating classes other than all Dirac measures $\mathcal D$ can be proved using Choquet theory, cf. \cite[Proposition 6.66]{shapiroLecturesStochasticProgramming2014} \cite[Lemma~3.1]{popescuSemidefiniteProgrammingApproach2005}, as well as in \cite{shapiroDualityTheoryConic2001,rogosinskiMomentsNonNegativeMass}.

\subsection{Compactness of the ambiguity set}
We now prove the compactness of the ambiguity set.
We use the mean map notation $\mmap : P \mapsto \mu_P$ to denote a map between the 
space of $\mathcal P$ equipped with MMD, and $\rkhs$ equipped with its norm.
Let us denote the image of a subset $\mathcal K$ of measures under \mmap by $\mmap(\mathcal{K}) := \{ \mu_P \mid P\in \mathcal{K}\} \subseteq \rkhs$.
If $\rkhs$ is universal, then MMD is a metric. 
By the definition of MMD, \mmap is an isometry (i.e., distance-preserving map) between $\mathcal{P}$ and $\rkhs$.
\begin{lemma}
        \label{thm:compact2}
    $\mmap (\mathcal P)$ is compact if \domain is compact.
\end{lemma}
\begin{proof}
If \domain is compact, by Prokhorov's theorem $\mathcal P$ is compact. Since $\mmap$ is an isometry, $\mmap(\mathcal P)$ is compact.
\end{proof}
It is straightforward to verify that $\mmap(\mathcal{P})$ is convex.
\begin{lemma}
    \label{thm:kreinmilman}
Let $C_P=\mathcal C\cap \mmap (\mathcal P)$. If \domain is compact, under Assumption~\ref{asm:compact}, $C_p$ is compact.
\end{lemma}
\begin{proof}
        By the Krein-Milman theorem, the convexity and compactness of $\mmap (\mathcal P)$ (proved in the previous lemma) imply that $\mmap (\mathcal P)$ is closed.
        By Assumption~\ref{asm:compact}, $\mathcal C$ is closed, which results in the closedness of $C_p$. Since $C_p$ is a closed subset of a compact set $\mmap(\mathcal P)$, it is compact.
\end{proof}

Recall that we denote the feasible set of probability measures, i.e., ambiguity set, for primal \kdro~\eqref{eq:kdro} by $ \pfeas = {\{P \colon \int \phi \ dP  = \mu, \mu\in\mathcal C, P\in \mathcal P\}}$. 
It is convex by straightforward verification. Let us derive the following compactness property of the ambiguity set.
\begin{lemma}
    If \domain is compact, under Assumption~\ref{asm:compact}, $\pfeas$ is compact.
    \label{thm:feas_is_compact}
\end{lemma}
\begin{proof}
    We first note $\pfeas=\mmap^{-1}(C_p)$ and \mmap is an isometric isomorphism (i.e., bijective isometry) between \pfeas and $C_p$.
    Then $\pfeas$ is compact since $C_p$ is compact.
\end{proof}